\documentclass[a4paper,11pt]{article}
\usepackage[latin1]{inputenc}
\usepackage[T1]{fontenc}
\usepackage{a4wide}
\usepackage{fontenc}
\usepackage{color}
\usepackage{amsmath}
\usepackage{graphicx}
\usepackage{array}
\usepackage{amsfonts}
\usepackage{amssymb}
\usepackage{euscript}
\usepackage{delarray}
\usepackage{latexsym}
\usepackage{enumerate}
\usepackage{float}

\begin{document}
\title{Recursive co-kriging model for Design of Computer experiments with multiple levels of fidelity with an application to hydrodynamic}
\author{Loic Le Gratiet \\ CEA, DAM, DIF, F-91297 Arpajon, France \\ loic.le-gratiet@cea.fr}

\maketitle

\section*{abstract}We consider in this paper the problem of building a fast-running approximation - also called surrogate model - of a complex computer code.  The co-kriging based surrogate model is a  promising tool to build such an approximation when the complex computer code can be run at different levels of accuracy. We present here an original  approach to perform a multi-fidelity co-kriging  model which is based on a recursive formulation.
We prove that the predictive mean and the variance   of the presented approach  are identical to the ones of the original co-kriging model proposed by [Kennedy, M.C. and O'Hagan, A., Biometrika, 87, pp 1-13, 2000]. However, our new approach allows   to obtain original  results.
First, closed form formulas for the  universal co-kriging predictive mean and variance are given.
Second,  a fast cross-validation procedure for the multi-fidelity co-kriging model is introduced. 
Finally, the proposed approach has a reduced    computational complexity compared to the  previous one.
The multi-fidelity model is successfully applied  to emulate a hydrodynamic simulator. \\

\paragraph{keywords}: uncertainty quantification, surrogate models, universal co-kriging, recursive model, fast cross-validation, multi-fidelity computer code.

\section{Introduction}

Computer codes are widely used in science and engineering  to describe physical phenomena. Advances in physics and computer science  lead to increased complexity for the simulators. 
As a consequence, to perform a sensitivity analysis or an optimization based on a complex computer code, a fast approximation of it - also called surrogate model - is built in order to avoid prohibitive computational cost.  A very popular method to build surrogate model is the Gaussian process regression, also named   kriging. It corresponds to  a particular class of surrogate models  which makes the assumption that  the response  of the complex code is  a realization of a Gaussian process. This method was originally introduced  in geostatistics  in \cite{Kri51}  and it was then proposed in the field of computer experiments  in \cite{SACKS89}. 
During the last decades, this method has become widely used and   investigated. The reader is referred to the books  \cite{S99}, \cite{S03} and \cite{R06} for more detail about it. 

Sometimes low-fidelity versions of the computer code are available. They may be less accurate but they are computationally cheap.
A question of interest is how to build a surrogate  model using data from simulations  of multiple levels of fidelity. Our  objective is hence to build a multi-fidelity surrogate model which is able to use the information obtained from the fast versions of the code. Such models have been presented in the literature \cite{Crai98}, \cite{KO00}, \cite{Hig04}, \cite{AF07},  \cite{QW07} and \cite{CG09}.

The first multi-fidelity model proposed in  \cite{Crai98} is based on a linear regression formulation. Then  this model is  improved in \cite{CG09} by using a Bayes linear formulation. The reader is referred to \cite{GW07} for further detail about the Bayes linear approach. The methods suggested in \cite{Crai98} and \cite{CG09} have the strength to be relatively  computationally clean  but as they are based on a linear regression formulation, they  could suffer from a lack of accuracy. Another  approach is to use an extension of  kriging  for multiple response models which is called co-kriging.  The idea is implemented in  \cite{KO00} which  presents  a co-kriging model based on an autoregressive relation between the different code levels. This method turns out to be very efficient and it has been applied and extended significantly. In particular,  the use of co-kriging for multi-fidelity optimization is presented in \cite{AF07}  and    a Bayesian formulation is proposed in \cite{QW07}.

The  strength of the co-kriging model is that it  gives   very good predictive models but it is often computationally expensive, especially when the number of simulations  is large. Furthermore,    large data set can generate  problems such as  ill-conditioned covariance matrices. These problems are known for kriging but they become even more difficult  for  co-kriging  since the total number of observations is the sum of the observations at all code levels.  

In this paper, we adopt a new approach for multi-fidelity surrogate modeling which uses a co-kriging model but with an original recursive formulation. In fact, our model is able to build a $s$-level co-kriging model by building $s$ independent krigings.
An important property of this model is that it provides   predictive mean and variance identical to the ones presented in \cite{KO00}.  However, our  approach significantly reduces the complexity of the model since it divides the whole set of simulations  into  groups of simulations  corresponding to the ones of each level. Therefore, we will have $s$ sub-matrices to invert which is less expensive and ill-conditioned than a large one and the estimation of the parameters can be performed separately (Section \ref{complexity}).

Furthermore, a strength  of our  approach is that it allows   to   extend classical results of kriging to the considered co-kriging model. 
The two original results presented in our paper are the following ones:
First, closed form expressions for the universal co-kriging predictive mean and variance are given (Section \ref{universal}).
Second,   the fast cross-validation method proposed in  \cite{Dub83} is extended to  the multi-fidelity co-kriging model (Section \ref{FastCV}). 
Finally, we illustrate these results in a complex hydrodynamic simulator (Section \ref{illustration}).

\section{Multi-fidelity  Gaussian process regression}\label{MufiGP}

In Subsection \ref{AR1model}, we briefly present the   approach to build multi-fidelity model suggested in \cite{KO00} that uses  a co-kriging model. In  Subsection \ref{recursivemodel},  we  detail our recursive approach to build  such a  model. The recursive formulation of the multi-fidelity model is the first novelty of this paper. We will see in the next sections that the new formulation allows us to find original  results about the co-kriging model and  to reduce its  computational complexity.  

\subsection{The classical autoregressive model}\label{AR1model}

Let us suppose that we have $s$ levels of code  $(z_t(x))_{t=1,\dots,s}$ sorted by increasing order of fidelity and  modeled by Gaussian processes $(Z_t(x))_{t=1,\dots,s}$, $x \in Q$. We hence consider that  $z_s(x)$ is  the most accurate and costly code that we want to surrogate and $(z_t(x))_{t=1,\dots,s-1}$ are cheaper versions of it with $z_1(x)$ the less accurate one.  We consider the following autoregressive model with  $ t = 2,\dots,s$:

\begin{equation}\label{eq1}
\left\{
\begin{array}{l}
Z_t(x) = \rho_{t-1}(x)Z_{t-1}(x)+\delta_t(x) \\
 Z_{t-1}(x) \perp \delta_t(x) \\
\rho_{t-1}(x) =g_{t-1}^T(x) \beta_{\rho_{t-1}} \\
\end{array}
\right.
\end{equation}
where:
\begin{equation}\label{eq2}
\delta_t(x) \sim \mathcal{GP} (f_t^T(x) \beta_t, \sigma_t^2 r_t(x,x'))
\end{equation}
and: 
\begin{equation}\label{eq3}
Z_1(x) \sim \mathcal{GP} (f_1^T(x) \beta_1, \sigma_1^2 r_1(x,x'))
\end{equation}
Here, $^T$ stands for the transpose, $\perp$ denotes the independence relationship, $\mathcal{GP}$ stands for Gaussian Process, $g_{t-1}(x)$ is a vector of $q_{t-1}$ regression functions, $f_t(x)$ is a vector of $p_t$ regression functions, $r_t(x,x')$ is a correlation function, $\beta_t$ is a $p_t$-dimensional vector, $\beta_{\rho_{t-1}}$ is a $q_{t-1}$-dimensional vector,  and $\sigma_t^2$ is  a real.
Since we  suppose that the responses are realizations of Gaussian processes, the multi-fidelity model can be built by conditioning   by the known responses of the codes at the different levels.

The previous model comes  from the article  \cite{KO00}. It is induced by the following assumption: $\forall x \in Q$,  if we know $Z_{t-1}(x)$, nothing more can be learned  about  $Z_t(x)$ from $Z_{t-1}(x')$ for  $x \neq x'$.

Let us consider $\mathcal{Z}^{(s)} = (\mathcal{Z}_1^T,\dots,\mathcal{Z}_s^T)^T$ the Gaussian vector containing the values of the   random processes $(Z_t(x))_{t=1,\dots,s}$ at the points in the experimental design sets $(D_t)_{t=1,\dots,s}$ with  $D_s \subseteq D_{s-1} \subseteq \dots \subseteq D_1$. We denote by  $z^{(s)} = (z_1^T,\dots,z_s^T)^T$ the  vector containing the values  of $(z_t(x))_{t=1,\dots,s}$ at the  points in $(D_t)_{t=1,\dots,s}$. 
The nested property of the experimental design sets is not necessary to build the model but it allows  for a simple  estimation of the model parameters. Since the codes are sorted in increasing order of fidelity it is not an unreasonable constraint for practical applications.
By denoting $\beta = (\beta_1^T,\dots,\beta_s^T)^T$ the trend parameters, $\beta_\rho = (\beta_{\rho_1}^T,\dots,\beta_{\rho_{s-1}}^T)^T$ the trend of the adjustment parameters and  $\sigma^2 = (\sigma_1^2,\dots,\sigma_s^2)$  the variance parameters, we have for any $x \in Q$:

\begin{displaymath} 
 [Z_s(x)|\mathcal{Z}^{(s)}=z^{(s)},\beta,\beta_\rho,\sigma^2] \sim \mathcal{N}\left(m_{Z_s}(x),s_{Z_s}^2(x)\right)
\end{displaymath}

where:
\begin{equation}\label{eq4}
m_{Z_s}(x) = h_s(x)^T\beta + t_s(x)^TV_s^{-1}(z^{(s)}-H_s\beta)
\end{equation}

and:
\begin{equation}\label{eq5}
s_{Z_s}^2(x) = v_{Z_s}^2(x) - t_s(x)^T V_s^{-1} t_s(x)
\end{equation}

The Gaussian process regression mean $m_{Z_s}(x)$ is the  predictive model of the highest fidelity response $z_s(x)$ which is built with the known responses of all code levels $z^{(s)}$. The variance $s_{Z_s}^2(x)$ represents the predictive mean squared error of the model.\\

The matrix $V_s$ is the covariance matrix of the Gaussian vector  $\mathcal{Z}^{(s)}$,  the vector $t_s(x)$ is the vector of covariances between $Z_s(x)$ and  $\mathcal{Z}^{(s)}$,  $H_s \beta$ is the mean of  $\mathcal{Z}^{(s)}$, $h_s(x)^T\beta$ is the mean of $Z_s(x)$ and  $v_{Z_s}^2(x)$ is the  variance of  $Z_s(x)$.  All these terms are built in terms of the experience vector at level $t$ (\ref{eq6}) and of  the covariance between $Z_t(x)$ and $Z_{t'}(x)$:

\begin{equation}\label{eq6}
h_t(x)^T = \left( \left(\prod_{i=1}^{t-1} {\rho_{i}(x) }\right)f_1^T(x) , \left(\prod_{i=2}^{t-1} {\rho_{i}(x) }\right)f_2^T(x) , \dots , \rho_{t-1}(x) f_{t-1}^T(x) , f_t^T(x) \right)
\end{equation}

\begin{equation}\label{eq7}
\mathrm{cov}(Z_t(x),Z_ {t'}(x')|\sigma^2,\beta,\beta_\rho) =
\left(\prod_{i=t'}^{t-1}\rho_i(x)\right) \mathrm{cov}(Z_{t'}(x),Z_{t'}(x')|\sigma^2,\beta,\beta_\rho)
\end{equation}
and:
\begin{equation}\label{eq7bis}
\mathrm{cov}(Z_t(x),Z_t(x')|\sigma^2,\beta,\beta_\rho) =
\sum_{j=1}^{t}{\sigma_{j}^2\left( \prod_{i=j}^{t-1} {\rho_{i}(x) \rho_{i}(x')} \right)r_j(x,x')}
\end{equation}

\paragraph{Remark.} The model (\ref{eq1}) is an extension of the model presented in  \cite{KO00} in which the adjustment parameters $\rho_t(x)_{t=2,\dots,s}$ do not depend on $x$. We show in a practical application (Section \ref{illustration})  that this extension is worthwhile.

\subsection{Recursive multi-fidelity model}\label{recursivemodel}

In this section, we present the  new multi-fidelity model which is based on a recursive formulation. Let us consider the following model for $t=2,\dots,s$ :

\begin{equation}\label{eq7ter}
\left\{
\begin{array}{l}
Z_t(x) = \rho_{t-1}(x)\tilde{Z}_{t-1}(x)+\delta_t(x) \\
\tilde{Z}_{t-1}(x)  \perp \delta_t(x) \\
\rho_{t-1}(x) = g_{t-1}^T(x) \beta_{\rho_{t-1}} \\
\end{array}
\right.
\end{equation}
where $\tilde{Z}_{t-1}(x)$ is a Gaussian process with distribution $[Z_{t-1}(x)|\mathcal{Z}^{(t-1)}=z^{(t-1)},\beta_{t-1},\beta_{\rho_{t-2}},\sigma_{t-1}^2]$, $\delta_t(x)$ is a Gaussian process with distribution (\ref{eq2})  and $D_s \subseteq D_{s-1} \subseteq \dots \subseteq D_1$. The unique difference with the previous model is that we express $Z_t(x)$ (the Gaussian process modeling the response at level $t$) as a function of the Gaussian process $Z_{t-1}(x)$ conditioned by the  values $z^{(t-1)}=(z_1,\dots,z_{t-1})$ at points in the experimental design sets $(D_i)_{i=1,\dots,t-1}$. 
As in the previous model, the nested property for the experimental design sets  is assumed because it allows for  efficient estimations of  the model parameters.
The Gaussian processes $(\delta_t(x))_{t=2,\dots,s}$ have the same definition as previously and we have for  $t=2,\dots,s$ and for $x \in Q$:
\begin{equation}\label{eq10}
\left[Z_{t}(x)|\mathcal{Z}^{(t)}=z^{(t)}, \beta_t,\beta_{\rho_{t-1}},\sigma_t^2 \right] \sim
\mathcal{N}
\left(
\mu_{Z_{t}}(x) ,
s_{Z_{t}}^2(x) 
\right)
\end{equation}
where:
\begin{equation}\label{eq11}
\mu_{Z_{t}}(x) = \rho_{t-1}(x) \mu_{Z_{t-1}}(x)  + f_t^T(x)\beta_t + 
r_t^T(x) R_t^{-1}\left(
z_t - \rho_{t-1}(D_t)\odot z_{t-1}(D_t) - F_t \beta_t
\right)
\end{equation}
and:
\begin{equation}\label{eq12}
\sigma_{Z_{t}}^2(x)  =  \rho_{t-1}^2(x)\sigma_{Z_{t-1}}^2(x) + \sigma_t^2\left(
1-r_t^T(x) R_t^{-1}r_t(x)
\right)
\end{equation}

The notation $\odot $ represents the element by element matrix product. $R_t$ is the correlation matrix  $R_t = (r_t(x,x'))_{x,x' \in D_t}$ and  $r_t^T(x)$ is the correlation vector $r_t^T(x) = (r_t(x,x'))_{x' \in D_t}$. We denote by $\rho_t(D_{t-1})$ the vector containing the values of  $\rho_t(x)$ for  $x\in D_{t-1}$, $z_t(D_{t-1})$ the vector containing the known values of $Z_t(x)$ at points in $D_{t-1}$  and  $F_t$ is the experience matrix  containing the values of $f_t(x)^T$ on $D_t$. \\

 The mean $\mu_{Z_{t}}(x)$ is the surrogate model of the response at level $t$, $1 \leq t \leq s$, taking into account the known values of the $t$ first levels of responses $(z_i)_{i=1,\dots,t}$ and the variance $\sigma_{Z_{t}}^2(x)$ represents the mean squared error of this model. The mean and the variance of the Gaussian process regression at level $t$ being expressed in function of the ones of level $t-1$, we   have a recursive multi-fidelity metamodel.   Furthermore, in this new formulation, it is clearly emphasized that the mean of the predictive distribution does not depend on the variance parameters $(\sigma_t^2)_{t=1,\dots,s}$. This is a classical result of kriging   which states that for covariance kernels of the form $k(x,x') = \sigma^2r(x,x')$, the mean of the kriging model is independent of $\sigma^2$.
Another  important strength of the recursive formulation is that  contrary to the formulation suggested in \cite{KO00}, once the multi-fidelity model is built,  it provides the surrogate models of   all the responses $(z_t(x))_{t=1\dots,s}$.

We have the following proposition.

\paragraph{Proposition 1}   Let us consider $s$ Gaussian processes $(Z_t(x))_{t=1,\dots,s}$ and $\mathcal{Z}^{(s)} = (\mathcal{Z}_t)_{t=1,\dots,s}$ the Gaussian vector  containing the values of $(Z_t(x))_{t=1,\dots,s}$ at points in  $(D_t)_{t=1,\dots,s}$ with  $D_s \subseteq D_{s-1} \subseteq \dots \subseteq D_1$. If we consider the mean and the variance (\ref{eq4}) and (\ref{eq5}) induced by the model (\ref{eq1}) when we condition the Gaussian process $Z_s(x)$ by the known values $z^{(s)}$ of $\mathcal{Z}^{(s)}$ and by the parameters $\beta$, $\beta_\rho$ and $\sigma^2$ and the mean and the variance (\ref{eq11}) and (\ref{eq12}) induced by the model (\ref{eq7ter})  when we condition $Z_s(x)$ by $z^{(s)}$ and by the parameters $\beta$, $\beta_\rho$ and $\sigma^2$, then, we have:
\begin{eqnarray*}
\mu_{Z_{s}}(x)  & =  & m_{Z_{s}}(x)  \\
\sigma_{Z_{s}}^2(x) & = & s_{Z_{s}}^2(x) 
\end{eqnarray*}

The proof of the proposition is given in Appendix \ref{A1}. It shows  that the model  presented in   \cite{KO00} and the recursive model (\ref{eq7ter}) have the same predictive Gaussian  distribution.  
Our objective in the next sections is to show that the new formulation (\ref{eq7ter}) has several advantages compared to the one of \cite{KO00}. First, its computational complexity is lower (Section \ref{complexity});  second, it provides closed form expressions for the universal co-kriging mean and variance contrarily to \cite{KO00} (Section \ref{universal}); third, it makes it possible to implement a fast cross-validation procedure (Section \ref{FastCV}).

\subsection{Complexity analysis}\label{complexity}

The computational cost is dominated by the inversion of the covariance matrices. In the original approach proposed in \cite{KO00} one has to invert the matrix $V_s$ of size $\sum_{i=1}^sn_i \times \sum_{i=1}^sn_i$.

Our recursive  formulation shows  that building a $s$-level co-kriging is equivalent to build $s$ independent krigings. This implies a reduction of the model complexity. Indeed, the inversion of 
$s$ matrices $(R_t)_{t={1,\dots,s}}$ of size $(n_t \times n_t)_{t={1,\dots,s}}$ where $n_t$ corresponds to the size of the vector $z_t$ at level $t = 1,\dots,s$  is less expensive than the inversion  of   the matrix $V_s$ of size $\sum_{i=1}^sn_i \times \sum_{i=1}^sn_i$. 
We also reduce the memory cost  since storing the $s$ matrices $(R_t)_{t={1,\dots,s}}$ requires less memory than storing the matrix $V_s$. Finally, we note that the model with this formulation is more interpretable since we  can deduce the impact of each level of response into the model error through  $(\sigma_{Z_{t}}^2(x))_{t={1,\dots,s}}$.

\subsection{Parameter estimation}\label{paramestim}

We present in this section a Bayesian estimation of the parameter  $\psi = (\beta, \beta_\rho, \sigma^2)$ focusing on conjugate and non-informative distributions  for the priors. This allows us to obtain closed form expressions for the estimations of the parameters. Furthermore, from the non-informative case, we can obtain the estimates given by a maximum likelihood method. The presented formulas can hence be used in a frequentist approach. We note that the recursive formulation  and the nested property of the experimental designs allow for separate the estimations of the parameters $(\beta_t,\beta_{\rho_{t-1}},\sigma^2_t)_{t=1,\dots,s}$ and $(\beta_1,\sigma_1^2)$.

We address two cases  in this section 
\begin{itemize}
\item Case (i): all the priors are informative
\item Case (ii):  all the priors are non-informative
\end{itemize}
It is of course  be possible to address the case of a mixture of informative and non-informative priors. For the non-informative case (ii), we use the ``Jeffreys priors''  \cite{J61}:
\begin{equation}
p(\beta_1|\sigma_1^2) \propto 1, \quad p(\sigma_1^2) \propto \frac{1}{\sigma_1^2}, \quad p(\beta_{\rho_{t-1}},\beta_t|z^{(t-1)},\sigma_t^2) \propto 1,  \quad p(\sigma_t^2|z^{(t-1)}) \propto  \frac{1}{\sigma_t^2}
\end{equation}
where $t=2,\dots,s$. For the informative case (i), we consider the following conjugate prior distributions:
\begin{displaymath}
[\beta_1|\sigma_1^2] \sim \mathcal{N}_{p_1}(b_1,\sigma_1^2V_1)
\end{displaymath}
\begin{displaymath}
 [\beta_{\rho_{t-1}},\beta_t| z^{(t-1)},\sigma_t^2] \sim \mathcal{N}_{q_{t-1}+p_t}\left(b_t = \left(\begin{array}{c}b^{\rho}_{t-1} \\ b^{\beta}_{t} \end{array} \right),\sigma_t^2V_t=\sigma_t^2 \left(\begin{array}{cc}V^{\rho}_{t-1} & 0 \\ 0 & V^{\beta}_{t} \\ \end{array} \right) \right)
\end{displaymath}
\begin{displaymath}
[\sigma_1^2] \sim \mathcal{IG}(\alpha_1,\gamma_1), \qquad [\sigma_t^2|z^{(t-1)}] \sim \mathcal{IG}(\alpha_t,\gamma_t)
\end{displaymath}
with $b_1$ a vector a size $p_1$, $b^{\rho}_{t-1}$ a vector of size $q_{t-1}$, $b^{\beta}_{t}$ a vector of size $p_t$, $V_1$ a $p_1\times p_1$ matrix, $V^{\rho}_{t-1}$ a $q_{t-1}\times q_{t-1}$ matrix, $V^{\beta}_{t}$ a $p_t \times p_t$ matrix,  $\alpha_1, \gamma_1, \alpha_t, \gamma_t > 0$ and $\mathcal{IG}$ stands for the inverse Gamma distribution. These informative priors allow the user to prescribe the prior means and   variances of all parameters. The choice of conjugate priors allows us  to have closed form expressions for the posterior distributions of the parameters. Indeed, we have:

\begin{equation}\label{eq8}
[\beta_1|z_1,\sigma_1^2] \sim \mathcal{N}_{p_1}(\Sigma_1 \nu_1, \Sigma_1) \qquad
[\beta_{\rho_{t-1}},\beta_t|z^{(t)}, \sigma_t^2] \sim \mathcal{N}_{q_{t-1}+q_t}(\Sigma_t \nu_t, \Sigma_t)
\end{equation}
where, for $t \geq 1$:
\begin{equation}\label{eq8ter}
\Sigma_t  = \left\{
\begin{array}{lr} 
\lbrack H_t^T \frac{R_t^{-1}}{\sigma_2^2}H_t+\frac{V_t^{-1}}{\sigma_2^2}\rbrack^{-1} &   \mathbf{(i)} \\ 
\lbrack H_t^T \frac{R_t^{-1}}{\sigma_2^2} H_t \rbrack^{-1} &\mathbf{(ii)}  \\
\end{array}
\right.
\nu_t = \left\{
\begin{array}{lr} 
\lbrack H_t^T \frac{R_t^{-1}}{\sigma_2^2}z_t+\frac{V_{t}^{-1}}{\sigma_2^2} b_t \rbrack &   \mathbf{(i)} \\
\lbrack H_t^T \frac{R_t^{-1}}{\sigma_2^2}z_t\rbrack & \mathbf{(ii)} \\
\end{array}
\right.
\end{equation}
with $H_1 = F_1$ and for $t > 1$,  $H_t=[G_{ t-1} \odot (z_{t-1}(D_t)\mathbf{1}^T_{q_{t-1}})\quad F_t] $ where $G_{t-1} $ is the experience matrix containing the values of $g_{t-1}(x)^T$ in $D_t$ and $\mathbf{1}^T_{q_{t-1}}$ is a $q_{t-1}$-vector of ones. Furthermore, we have for $t \geq 1$:
\begin{equation}\label{eq9}
[\sigma_t^2|z^{(t)}] \sim \mathcal{IG}(a_t, \frac{Q_t}{2})
\end{equation}
where:
\begin{displaymath}
Q_t = \left\{
\begin{array}{lr}
\gamma_t+(b_t-\hat{\lambda}_t)^T(V_t+[H_t^TR_t^{-1}H_t]^{-1})^{-1}(b_t-\hat{\lambda}_t)+\hat{Q}_t & \mathbf{(i)}  \\
\hat{Q}_t  & \mathbf{(ii)} \\
\end{array}
\right.
\end{displaymath}
with $\hat{Q}_t = (z_t-H_t \hat{\lambda}_t) ^T R_t^{-1} (z_t-H_t \hat{\lambda}_t)  $ ,  $ \hat{\lambda}_t=(H_t^TR_t^{-1}H_tF)^{-1}H_t^TR_t^{-1}z_t$ and :
\begin{displaymath}
a_t = \left\{
\begin{array}{lr}
\frac{n_t}{2}+\alpha_t &  \mathbf{(i)}\\
\frac{n_t-p_t-q_{t-1}}{2} &  \mathbf{(ii)}\\
\end{array}
\right.
\end{displaymath}
 with the convention $q_0 = 0$.

We highlight that the maximum likelihood   estimators for the parameters $\beta_1$ and $(\beta_{\rho_{t-1}}, \beta_t)$ are given by the means of the posterior distributions of the Bayesian estimations  in the non-informative case. Furthermore,  the  restricted maximum likelihood estimate of the variance parameter $\sigma_t^2$ can also be deduced from the  posterior distribution  of the Bayesian estimation  in the non-informative case and is given by 
$
\hat{\sigma}_{t,\mathrm{EML}}^2 = \frac{Q_t}{2a_t}
$.
The restricted maximum likelihood estimation  is a method which allows   to reduce the bias of the maximum likelihood estimation  \cite{PT71}.

\section{Universal co-kriging model}\label{universal}

We can see in equation (\ref{eq10}) that the predictive distribution of $Z_s(x)$ is conditioned by the observations $z^{(s)}$ and the parameters $\beta$, $\beta_\rho$ and $\sigma^2$.  The objective of a Bayesian prediction is to integrate the uncertainty due to the parameter estimations into   the predictive distribution. Indeed, in the previous subsection, we have expressed the posterior distributions of the variance parameters  $(\sigma_t^2)_{t=1,\dots,s}$ conditionally to the observations and the posterior distributions of the trend parameters $\beta_1$ and $(\beta_{\rho_{t-1}},\beta_t)_{t=2,\dots,s}$ conditionally to the observations and the variance parameters. Thus, using the Bayes formula, we can easily obtain a predictive distribution only conditioned by  the observations by integrating into it the posterior distributions of the parameters. 

As a result of this integration, the  predictive distribution is   not Gaussian. In particular, we cannot have a closed form expression for the predictive distribution. However, it is possible  to  obtain closed form expressions for  the posterior mean $\mathbb{E}[Z_s(x)|\mathcal{Z}^{(s)} = z^{(s)}]$ and   variance $\mathrm{Var}(Z_s(x)|\mathcal{Z}^{(s)} = z^{(s)})$. 

The following proposition giving the closed form expressions of the posterior  mean and   variance of the predictive distribution only conditioned by the observations is a novelty. The proof of this proposition is based on the recursive formulation which emphasizes the strength  of this new approach. Indeed,  it does not seem possible to obtain  this   result by considering directly the model suggested  in \cite{KO00}.

\paragraph{Proposition 2}  Let us consider $s$ Gaussian processes $(Z_t(x))_{t=1,\dots,s}$ and $\mathcal{Z}^{(s)} = (\mathcal{Z}_t)_{t=1,\dots,s}$ the Gaussian vector  containing the values of $(Z_t(x))_{t=1,\dots,s}$ at points in  $(D_t)_{t=1,\dots,s}$ with  $D_s \subseteq D_{s-1} \subseteq \dots \subseteq D_1$. If we consider the conditional predictive distribution in equation (\ref{eq10}) and the posterior distributions of the parameters given in equations (\ref{eq8}) and (\ref{eq9}), then we have for $t=1,\dots,s$:
\begin{equation}\label{13}
 \mathbb{E}[Z_t(x)|\mathcal{Z}^{(t)} = z^{(t)}] = 
h_t^T(x)
\Sigma_t \nu_t + r_t^T(x)R_t^{-1} \left( z_t - H_t \Sigma_t \nu_t  \right)
\end{equation}
with $h_1^T = f_1^T$, $H_1 = F_1$ and for $t > 1$,  $h_t^T(x) = \left(
\begin{array}{lr}
g_{t-1}(x)^T\mathbb{E}[Z_{t-1}(x)|\mathcal{Z}_{t-1}  = z_{t-1} ]   &  f_{t}^T(x)\\
\end{array}
\right)$ and $H_t=[G_{t-1} \odot (z_{t-1}(D_t)\mathbf{1}^T_{q_{t-1}})\quad F_t] $.
Furthermore, we have:
\begin{equation}\label{19}
\begin{array}{lll}
\mathrm{Var}(Z_t(x)|\mathcal{Z}^{(t)} = z^{(t)}) &  = &  \hat{\rho}_{t-1}^2(x) \mathrm{Var}(Z_{t-1}(x )|\mathcal{Z}^{(t-1)} = z^{(t-1)} )+\frac{Q_t}{2\left( a_t-1 \right)}\left(1-r_t^T(x)R_t^{-1}r_t^T(x)\right) \\
 & &+ \left(
h_t^T - r_t^T(x)R_t^{-1}H_t\right) \Sigma_t \left(
h_t^T - r_t^T(x)R_t^{-1}H_t\right)^T \\
\end{array}
\end{equation}
with  $\hat{\rho}_{t-1}(x) = [\Sigma_t \nu_t]_{1,\dots,q_{t-1}}$.

The proof of Proposition 2 is given in Appendix \ref{A2}. We note that, in the mean of the predictive distribution, the parameters have been  replaced by their posterior means. Furthermore, in the variance of the predictive distribution, the variance parameter  has been replaced by its  posterior mean  and the  term $\left(
h_t^T - r_t^T(x)R_t^{-1}H_t\right) \Sigma_t \left(
h_t^T - r_t^T(x)R_t^{-1}H_t\right)^T$   has been added.
 It represents the uncertainty due to the estimation of the regression parameters (including the adjustment coefficient).
We  call these formulas the universal co-kriging equations due to their similarities with the well-known universal kriging equations (they are identical for $s=1$).

\section{Fast cross-validation for   co-kriging surrogate models}\label{FastCV}

The idea of a cross-validation procedure is to split the experimental design  set into two disjoint  sets, one is used for training and the other one  is used to monitor the performance of the surrogate  model. The  idea is that the performance on the test set can be used as a proxy  for the generalization error. A particular case of this method is the Leave-One-Out Cross-Validation (noted  LOO-CV) where $n$ test sets are obtained by removing one observation at a time. This procedure can be  time-consuming for a kriging model  but it is shown in  \cite{Dub83}, \cite{R06} and \cite{ZW09}    that there are computational shortcuts.  Our recursive formulation allows to extend these ideas to co-kriging models (which is not possible with the original formulation in \cite{KO00}). Furthermore, the cross-validation equations proposed in this section extend the previous ones even for $s=1$  (i.e. the classical kriging model) since they do not suppose that the regression and the variance coefficients are known. Therefore, those parameters are re-estimated for each training set. We  note that the re-estimation of the variance coefficient is a novelty which is important since fixing this parameter can lead to huge  errors for the estimation  of the cross-validation predictive variance when the number of observations is small or when the number of points in the test set is important.

If we denote by $\xi_s$ the  set of indices   of  $n_{test}$ points in $D_s$ constituting the test set $D_{test}$ and $\xi_t$,  $1 \leq t < s$,  the corresponding  set of indices  in $D_t$ - indeed, we have $D_s \subset D_{s-1} \subset \dots \subset D_1$, therefore $D_{test} \subset D_t$. The nested experimental design assumption implies that, in the cross-validation procedure, if we remove a set of  points  from $D_s$ we  can also remove it from $D_t$, $1 \leq t \leq  s$.

The following proposition gives the vectors of the cross-validation predictive errors  and variances  at points in the  test set $D_{test}$ when we remove them  from the $t$ highest levels of code. In the proposition, we consider that we are in the non-informative case for the parameter estimation (see Section \ref{paramestim}) but it can be easily extended to the informative case presented in Section \ref{paramestim}. We note that this result presented for the first time to a multi-fidelity co-kriging model can be obtained thanks to the recursive formulation.

\paragraph{Notations:} If $\xi$ is a set of indices, then $A_{[\xi,\xi]}$ is the sub-matrix of elements $\xi \times \xi$ of $A$, $a_{[\xi]}$ is the sub-vector of elements $\xi$ of  $a$, $B_{[-\xi]}$ represents the matrix $B$ in which we remove   the rows of index $\xi$, $C_{[-\xi,-\xi]}$ is the sub-matrix of $C$ in which we remove the rows and columns of index $\xi$ and $C_{[-\xi,\xi]}$ is the sub-matrix of $C$ in which we remove the rows of index $\xi$  and keep the columns of index $\xi$.

\paragraph{Proposition 3} 
Let us consider $s$ Gaussian processes $(Z_t(x))_{t=1,\dots,s}$ and $\mathcal{Z}^{(s)} = (\mathcal{Z}_t)_{t=1,\dots,s}$ the Gaussian vector  containing the values of $(Z_t(x))_{t=1,\dots,s}$ at points in  $(D_t)_{t=1,\dots,s}$ with  $D_s \subseteq D_{s-1} \subseteq \dots \subseteq D_1$. We note $D_{test}$ a set made with the points of index $\xi_s$ of $D_s$ and  $\xi_t$ the corresponding points in $D_t$ with $1 \leq t \leq s$. Then, if we note $\varepsilon_{Z_s,\xi_s}$  the errors (i.e. real values minus predicted values) of the cross-validation procedure when we remove the points of $D_{test}$ from the $t$ highest levels of code, we have:
\begin{equation}\label{CVerr}
\begin{array}{lll}
\left(\varepsilon_{Z_s,\xi_s} - \rho_{s-1}(D_{test}) \odot \varepsilon_{Z_{s-1},\xi_{s-1}} \right) \left[ R_s^{-1} \right]_{[\xi_s,\xi_s]} & = & 
\left[ R_s^{-1}\left(
z_s - H_s \lambda_{s,-\xi_s}
\right) \right]_{[\xi_s ]} \\
\end{array}
\end{equation}
with $\varepsilon_{Z_u,\xi_u} = 0$ when  $u < t$, $\lambda_{s,-\xi_s} \left([H_s^T ]_{[-\xi_s]}K_s [H_s]_{[-\xi_s]} \right)= [H_s^T]_{[-\xi_s]} K_s z_s(D_s \setminus D_{test})$ and:
\begin{equation}\label{KtCV}
K_s = \left[ R_s^{-1} \right]_{[-\xi_s,-\xi_s]}-\left[ R_s^{-1} \right]_{[-\xi_s,\xi_s]}\left( \left[ R_s^{-1} \right]_{[\xi_s,\xi_s]} \right)^{-1}\left[ R_s^{-1} \right]_{[\xi_s,-\xi_s]}
\end{equation}
Furthermore, if we note $\sigma_{Z_s,\xi_s}^2$ the variances of the corresponding cross-validation procedure, we have:
\begin{equation}\label{CVvar}
\sigma_{Z_s,\xi_s}^2 = \rho_{s-1}^2(D_{\mathrm{test}}) \odot \sigma_{Z_{s-1},\xi_{s-1}}^2 + \sigma_{s,-\xi_s}^2\mathrm{diag}\left(
\left( \left[ R_s^{-1} \right]_{[\xi_s,\xi_s]}\right)^{-1}
\right)+\mathcal{V}_s
\end{equation}
with:
\begin{equation}\label{sigCV}
 \sigma_{s,-\xi_s}^2 = \frac{\left(z_s(D_s \setminus D_{\mathrm{test}}) - [H_s]_{[-\xi_s]}\lambda_{s,-\xi_s}  \right)^T K_s \left(z_s(D_s \setminus D_{\mathrm{test}}) - [H_s]_{[-\xi_s]}\lambda_{s,-\xi_s}  \right)  }{n_s-p_s-q_{s-1}-n_{train}}
\end{equation}
where $ \sigma_{u,-\xi_u}^2 = 0$ when $u < t$, $n_{train}$ is the length of the index vector $\xi_s$,  $H_s=[G_{s-1} \odot (z_{s-1}(D_s)\mathbf{1}^T_{q_{s-1}})\quad F_s] $ and:
\begin{equation}\label{CVUnivCok}
\mathcal{V}_s =
\mathcal{U}_s^T
\left([H_s^T ]_{[-\xi_s]}K_s [H_s]_{[-\xi_s]} \right)^{-1}
\mathcal{U}_s
\end{equation}
with $\mathcal{U}_s =   \left( \left( [R_s^{-1}]_{[\xi_s,\xi_s]}\right)^{-1}
\left[ R_s^{-1}H_s \right]_{[\xi_s ]} 
\right)$.

We note that  these equations are also valid  when $s=1$, i.e. for kriging model. We hence have closed form expressions  for the equations of a $k$-fold cross-validation with a re-estimation of the regression and variance parameters. These expressions can be deduced  from the universal co-kriging equations. The complexity of this procedure is essentially determined by the inversion of the matrices $\left(\left[ R_u^{-1} \right]_{[\xi_u,\xi_u]}\right)_{u=t,\dots,s}$  of size $n_{test} \times n_{test}$.  Furthermore, if we suppose the  parameters of variance and/or trend as known, we do not have to compute $\sigma_{t,-\xi_t}^2$ and/or  $\lambda_{t,-\xi_t}$ (they are fixed to their estimated value, i.e. $\sigma_{t,-\xi_t}^2 = \frac{Q_t}{2(a_t-1)}$ and $\lambda_{t,-\xi_t} = \Sigma_t \nu_t$, see Section \ref{paramestim})   which reduces substantially the complexity of the method. These equations generalize those of \cite{Dub83} and \cite{ZW09} where the variance $\sigma_{t,-\xi_t}^2$  is supposed to be  known. 
Finally, the term $\mathcal{V}_s$ is the additive  term due to the parameter estimations in the universal co-kriging. Therefore, if the trend parameters are supposed to be known, this term is equal to 0.
The proof of Proposition 3 is given in Appendix \ref{A3}.

\paragraph{Remark:} We must recognize that our closed form cross-validation formulas do not allow for the re-estimation of the hyper-parameters of the correlation  functions. However, as discussed in Subsection \ref{hyperparamestim}, Proposition 3  is useful even in that case to reduce the computational complexity of the cross-validation procedure.


\section{Illustration: hydrodynamic simulator}\label{illustration}

In this section we apply our  co-kriging method    to the hydrodynamic code   ``MELTEM''. 
The aim of the study is to build a prediction as accurate as possible using only a  few runs of the complex code and to assess the uncertainty of this prediction. In particular, we show the efficiency of the co-kriging model compared to the kriging one. We also illustrate the difference between simple and universal co-kriging and the results of the LOO-CV procedure. 
These illustrations are made possible and easy by the closed form formulas for the predictive mean and variance for universal co-kriging and by the fast cross-validation procedure described in Section \ref{FastCV} and \ref{universal} respectively.
Finally, we show that considering  an adjustment coefficient depending on $x$ can be worthwhile.

The  code MELTEM  simulates a second-order turbulence model for gaseous mixtures induced by Richtmyer-Meshkov instability \cite{MELTEM}. Two input parameters $x_1$ and $x_2$ are considered. They are phenomenological coefficients used in the equations of the energy of dissipation of the turbulent flow. These two coefficients vary in  the region $[0.5,1.5]\times[1.5,2.3]$. The considered code outputs, called $eps$ and $L_c$, are respectively  the dissipation factor and the mixture characteristic length. The simulator is a finite-elements code which can be run at  $s=2$ levels of accuracy by altering the finite-elements mesh. The simple code $z_1(.)$, using a coarse mesh, takes 15 seconds to produce an output whereas the complex code $z_2(.)$, using a fine mesh, takes 8 minutes. We use $5$ runs for the complex code $z_2(x)$ and $25$ runs for the cheap code $z_1(x)$. This represents 8 minutes on a hexa-core processor, which is our constraint for an operational use. Then, we build an additional set of $175$ points to test the accuracy of the models.  We note that  no prior information is available: we are hence  in the non-informative case.

\subsection{Estimation of the hyper-parameters}\label{hyperparamestim}

In the previous sections, we have considered  the correlation kernels $(r_t(x,x'))_{t=1,\dots,s}$ as known. In practical applications, we choose these kernels in a parameterized family of correlation kernels. Therefore, we consider kernels such that $r_t(x,x') = r_t(x,x';\phi_t)$. For $t=1,\dots,s$  the hyper-parameter $\phi_t$ can be estimated by maximizing the concentrated restricted log-likelihood \cite{S03} with respect to $\phi_t$:
\begin{equation}\label{reml}
\mathrm{log}\left( \left| \mathrm{det}\left(R_t\right) \right| \right)+ \left(n_t - p_t - q_{t-1}\right) \mathrm{log}\left( \sigma^2_{t, \mathrm{EML}} \right)
\end{equation}
with the convention $q_0=0$ and $ \sigma^2_{t, \mathrm{EML}}$ is the restricted likelihood estimate of the variance $\sigma_t^2$ (see Section \ref{paramestim}). This minimization problem has to be solved numerically.

 It is a common choice to estimate the hyper-parameters  by maximum likelihood \cite{S03}. 
It is also possible to estimate the hyper-parameters $(\phi_t)_{t=1,\dots,s}$ by minimizing a loss function of a Leave-One-Out Cross-Validation procedure. Usually, the complexity of this procedure is $\mathcal{O}\left( \left( \sum_{i=1}^s n_i\right)^4\right)$. Nonetheless, thanks to Proposition 3, it is reduced to $\mathcal{O}\left( \sum_{i=1}^s n_i^3 \right)$ since  it is essentially determined  by the inversions of the $s$ matrices  $(R_t )_{t=1,\dots,s}$.Therefore, the complexity for the estimation of $(\phi_t)_{t=1,\dots,s}$ is substantially reduced. Furthermore, the recursive formulation of the problem allows us to estimate the parameters  $(\phi_t)_{t=1,\dots,s}$ one at a time by starting with $\phi_1$ and estimating $\phi_t$, $t=2,\dots,s$ recursively.

\subsection{Comparison between  kriging and multi-fidelity co-kriging}

Before considering the real case study, we propose in this section a comparison between the kriging and co-kriging models when the number of runs $n_2$ for the complex code  varies such that $n_2 = 5, 10, 15, 20, 25$. For the co-kriging model, we consider    $n_1 = 25$ runs for the cheap code.
In this section, we focus on the output  $eps$. 

To perform the comparison,  we generate randomly 500 experimental design sets $(D_{2,i},D_{1,i})_{i=1,\dots,500}$ such that $D_{2,i} \subset D_{1,i}$, $i=1,\dots,500$, $D_{1,i}$ has $n_1$ points and $D_{2,i}$ has $n_2$ points.

We use for both  kriging and co-kriging models  a Matern$\frac{5}{2}$ covariance kernel and we consider $\rho$, $\beta_1$ and $\beta_2$ as constant.  The accuracies of the two models are 
evaluated on the test set composed of 175 observations. From them, the Root Mean Squared Error (RMSE) is computed: $\mathrm{RMSE} = \left(\frac{1}{175}\sum_{i=1}^{175} (\mu_{Z_2}(x_i^{\mathrm{test}}) - z_2(x_i^{\mathrm{test}}))^2\right)^{1/2}$.

\begin{figure}[h]
\begin{center}
\includegraphics[width = 8cm]{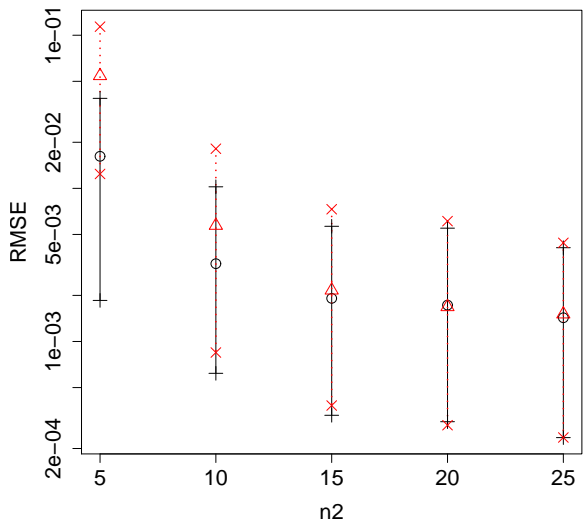}
\caption{Comparison between  kriging and co-kriging with $n_1 = 25$ runs for the cheap code (500 nested design sets have been randomly  generated for each $n_2$). The circles represent the averaged RMSE of the co-kriging, the triangles represent the averaged RMSE of the kriging, the crosses  represent the quantiles of probability 5\% and 95\% for the co-kriging RMSE and the  times signs represent the quantiles of probability 5\% and 95\% of the kriging RMSE. Co-kriging predictions are better than the ordinary kriging ones for small $n_2$ and they converge to the same accuracy when $n_2$ tends to $n_1 = 25$.}
\label{fig:comparison}
\end{center}
\end{figure}

Figure \ref{fig:comparison} gives the mean and the quantiles of probability 5\% and 95\% of the RMSE computed from the 500 sets $(D_{2,i},D_{1,i})_{i=1,\dots,500}$ when the number of runs for the expensive code $n_2$ varies.
In  Figure \ref{fig:comparison}, we can see that the errors converge  to the same value when $n_2$ tends to $n_1$. Indeed, due to the Markov property given in Section \ref{AR1model}, when $D_2 = D_1$, only the observations $z_2$ are taken  into account. Furthermore, we  can see that for small values of $n_2$, it is worth considering the co-kriging model since its accuracy  is significantly better than the one of the kriging model.

\subsection{Nested space filling design}\label{design}

As presented in Section \ref{MufiGP} we consider nested experimental design sets: $\forall t=2,\dots,s \quad D_t \subseteq D_{t-1}$. Therefore, we have to adopt particular design strategies to uniformly spread the inputs for all $D_t$. A strategy based on Orthogonal array-based Latin hypercube for nested space-filling designs is proposed by \cite{Q09a}.\\
We consider here another strategy for space-filling design, described in the following algorithm, which is very simple and not time-consuming. The number of points $n_t$ for each design $D_t$ is prescribed by the user, as well as the experimental design method applied to determine the coarsest grid $D_s$ used for the most expensive code $z_s$ (see \cite{FAN06} for a review of different methods). 

ALGORITHM 
\begin{itemize}
\item[] build $D_s = \{ x_j^{(s)} \}_{j=1,\dots,n_s} $ with the experimental design method prescribed by the user.
\item[] \textbf{for} $t$ = $s$ to $2$ \textbf{do}:
\begin{itemize}
\item[] build design $\tilde{D}_{t-1}$ with the experimental design method prescribed by the user.
\item[] \textbf{for} $i$ = $1$ to $n_t$ \textbf{do}:
\begin{itemize}
\item[] find $\tilde{x}_j^{(t-1)} \in \tilde{D}_{t-1}$ the closest point from $x_i^{(t)} \in D_t$ where $j \in [1, n_{t-1}]$.
\item[] remove $\tilde{x}_j^{(t-1)}$ from $\tilde{D}_{t-1}$.
\end{itemize}
\item[] \textbf{end for}
\item[] $D_{t-1} = \tilde{D}_{t-1} \cup D_t$.
\end{itemize}
\item[] \textbf{end for}
\end{itemize}
This strategy allows us to use any space-filling design method and it conserves the initial structure of the experimental design $D_s$ of the most accurate code, contrarily to a strategy based on selection of subsets of an experimental design for the less accurate code as presented by 
\cite{KO00}  and \cite{AF07}. We hence can ensure that $D_s$ has excellent space-filling properties. Moreover, the experimental design $D_{t-1}$ being equal to $\tilde{D}_{t-1} \cup D_t$, this method ensure the nested property.

\begin{figure}[h]
\begin{center}
\includegraphics[width = 8cm]{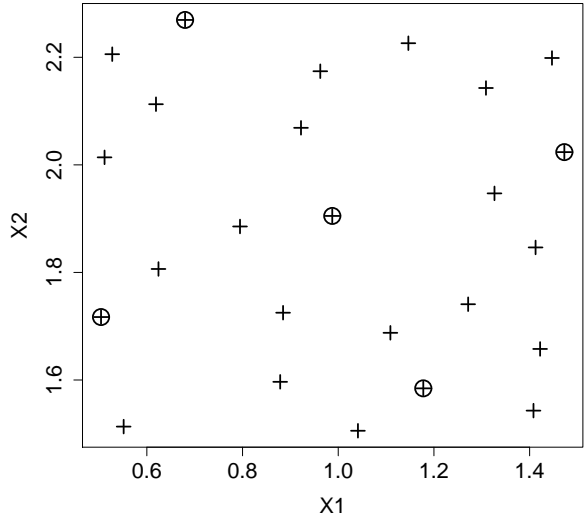}
\caption{Nested experimental design sets for the hydrodynamic application. The crosses represent the $n_1 = 25$ points of the experimental design set $D_1$ of the cheap code and the circles represent the $n_2 = 5$ points of the experimental design set $D_2$ of the expensive code.}
\label{fig:plan}
\end{center}
\end{figure}

In the presented application, we consider  $n_2 = 5$ points for the expensive code $z_2(x)$ and $n_1 = 25$ points for the cheap one $z_1(x)$. We apply the previous algorithm to build $D_2$ and $D_1$ such that $D_2 \subset D_1$. For the experimental design set $D_2$, we use a Latin-Hypercube-Sampling   \cite{St87} optimized with respect to the S-optimality criterion which  maximizes the mean distance from each design point to all the other points \cite{Sto05}. Furthermore, the set $D_1$ is built using a  maximum entropy design  \cite{Sh87} optimized with  the  Fedorov-Mitchell exchange algorithm \cite{Cu91}.  These algorithms are implemented in the library R lhs.
The obtained nested designs are shown in Figure \ref{fig:plan}.

\subsection{Multi-fidelity surrogate model for the dissipation factor $eps$}

We build here a co-kriging model for the dissipation factor $eps$. The obtained model is compared to a kriging one. This first example is used to illustrate the efficiency of the co-kriging method compared to  the kriging. It will also allow  us to highlight the difference between the simple and the universal co-kriging.

We use the experimental design sets presented in Section \ref{design}. To validate and compare our models, the 175 simulations of the complex code  uniformly spread on $[0.5,1.5]\times[1.5,2.3]$  are used. To build the different correlation matrices, we consider a tensorised matern-$\frac{5}{2}$ kernel (see \cite{R06}):
\begin{equation}\label{mat2D}
r(x,x';\theta_t) = r_{1d}(x_1,x_1';\theta_{t,1})  r_{1d}(x_2,x_2';\theta_{t,2}) 
\end{equation}
with $x=(x_1,x_2) \in [0.5,1.5]\times[1.5,2.3]$, $\theta_{t,1}, \theta_{t,2} \in (0,+\infty)$ and:
\begin{equation}
r_{1d}(x_i,x_i';\theta_{t,i}) = \left( 1 + \sqrt{5}\frac{|x_i-x_i'|}{\theta_{t,i}}+\frac{5}{3}\frac{(x_i-x_i')^2}{\theta_{t,i}^2}\right)\mathrm{exp}\left(-\sqrt{5}\frac{|x_i-x_i'|}{\theta_{t,i}}\right)
\end{equation}
Then, we consider $g_1(x) = 1$, $f_2(x) = 1$, $f_1(x) = 1$ (see Section \ref{AR1model} and \ref{recursivemodel}) and, using the concentrated maximum likelihood (see Section \ref{hyperparamestim}), we have the following estimations for the correlation hyper-parameters: $\hat{\theta}_1 = (0.69, 1.20) $ and  $\hat{\theta}_2 = (0.27, 1.37)$.

According to the values of the hyper-parameter estimates, the co-kriging model  is    smooth since the correlation lengths  are of the same order as the size of the input parameter space. Furthermore, the estimated correlation between the two codes is $82.64\%$, which shows that  the amount of information contained in the cheap code is  substantial.

Table \ref{tab3b} presents the results  of the parameter estimations (see Section \ref{paramestim}).

\begin{table}[H]
\begin{center}
\begin{tabular}{|c|c|c|}
\hline
Trend coefficient & $\Sigma_t \nu_t$ & $\Sigma_t/\sigma_t^2$ \\
\hline
$\beta_1$ & $8.84$ & $ 0.48$ \\
\hline
$\left( \begin{array}{c} \beta_{\rho_1} \\ \beta_2 \\ \end{array} \right)$ & $ \left( \begin{array}{c} 0.92  \\ 0.74  \\ \end{array} \right)$ & $ \left( \begin{array}{cc} 1.98   & -18.13  \\  -18.13 & 165.82 \\  \end{array} \right) $\\
\hline
\hline
Variance coefficient & $Q_t$ & $2 \alpha_t$ \\
\hline
$\sigma_1^2$ & $6.98$ & $24$ \\
\hline
$\sigma_2^2$ & $0.06$ & $3$ \\
\hline
\end{tabular}
\end{center}
\caption{Application: hydrodynamic simulator. Parameter estimation results  for the response $eps$ (see equations (\ref{eq8}) and (\ref{eq9})).}
\label{tab3b}
\end{table}

We see in Table \ref{tab3b} that the correlation between $\beta_{\rho_1}$ and $\beta_2$ is important which highlights the importance of taking into account the correlation between these two coefficients for the parameter estimation. We also see that the adjustment parameter $\beta_{\rho_1}$ is close to 1,  which means that the two codes are highly correlated. 

\begin{figure}[h]
\begin{center}
\includegraphics[width =7cm]{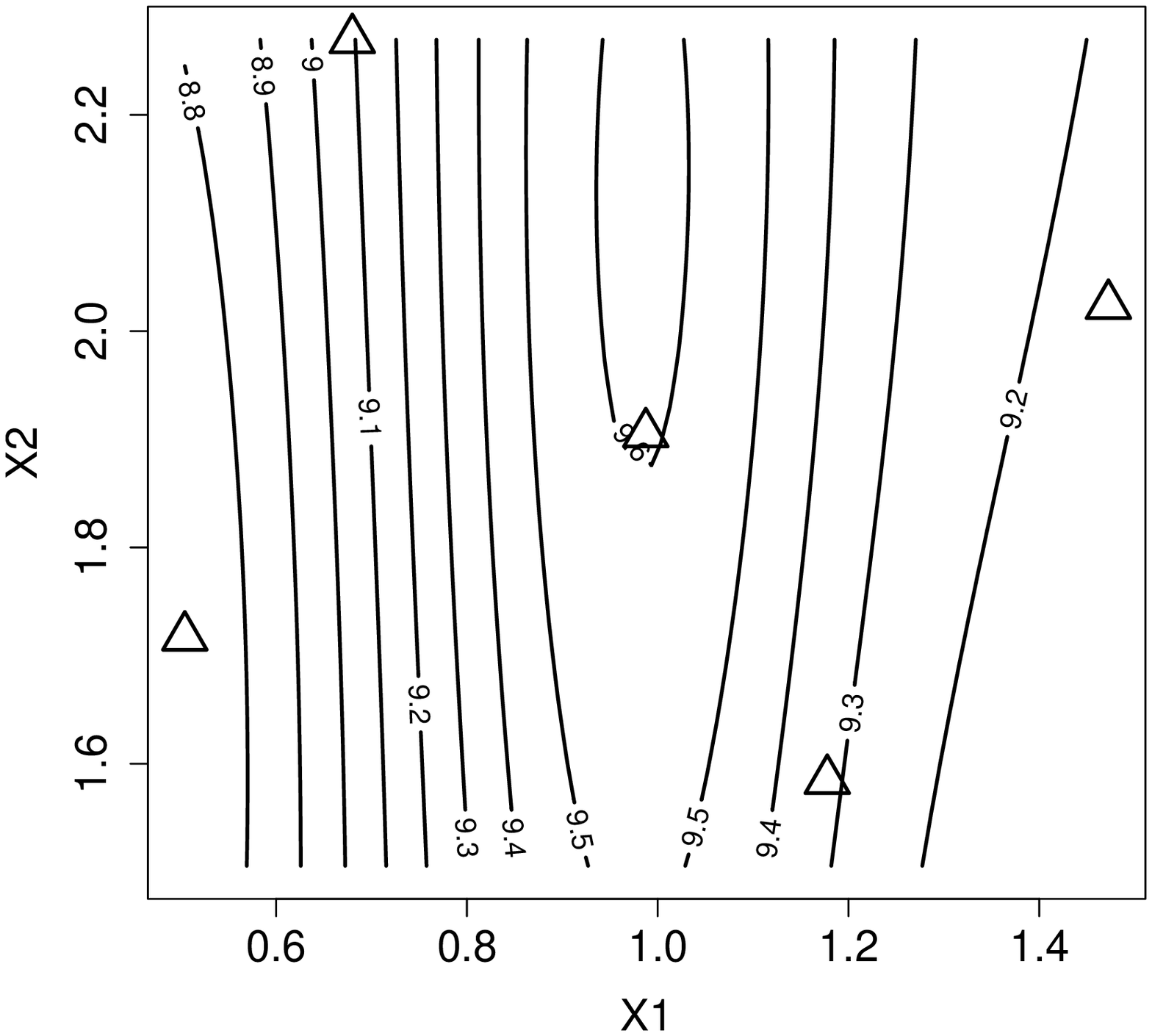}
\includegraphics[width = 7cm]{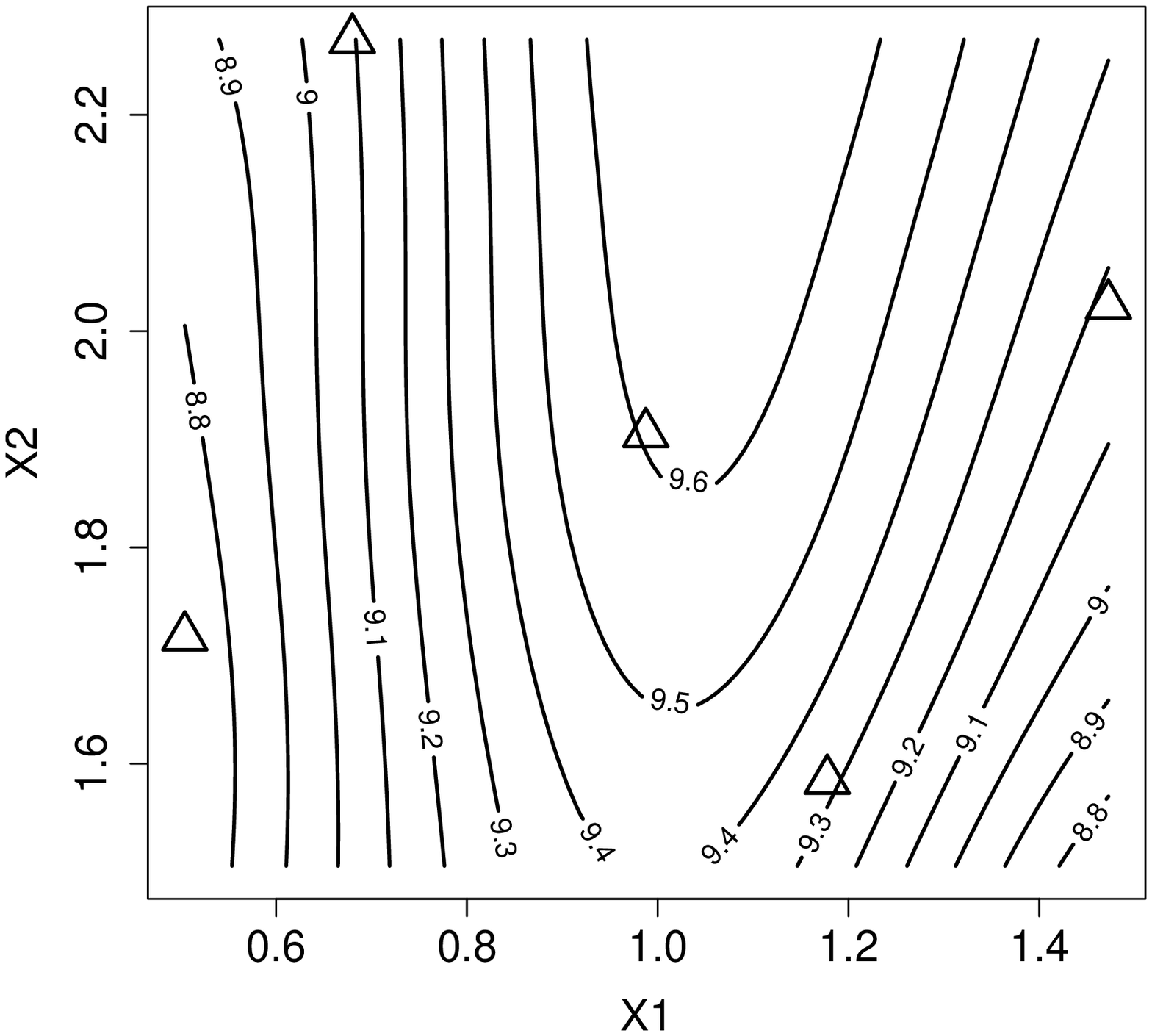}
\caption{Contour plot of the kriging mean (left picture) and the co-kriging mean (right picture). The triangles represent the $n_2 = 5$ points of the experimental design set  of the expensive code.}
\label{fig:contoureps}
\end{center}
\end{figure}

 Figure \ref{fig:contoureps} illustrates the contour plot of the kriging and co-kriging means, we can see significant  differences between the two surrogate models.

Table \ref{tab4} compares the prediction accuracy of the co-kriging and the kriging models. The different coefficients are estimated with the 175 responses of the complex code on the test set:
\begin{itemize}
\item MaxAE:  Maximal absolute value of the observed error.
\item RMSE : Root mean squared value of the observed error.
\item $Q_2 = 1 - ||\mu_{Z_2}(D_\mathrm{test})-z_2(D_\mathrm{test})||^2/||\mu_{Z_2}(D_\mathrm{test})-\bar{z}_2||^2$, with $\bar{z}_2 = (\sum_{i=1}^{n_2}z_2(x_i^\mathrm{test}))/n_2$.
\item RIMSE : Root of the average value of the kriging or co-kriging variance.
\end{itemize}
\begin{table}[H]
\begin{center}
\begin{tabular}{ccccc}
\hline
& $Q_2$ & RMSE & MaxAE  & RIMSE. \\
kriging & 75.83\% & 0.133 & 0.49  &  0.110 \\
co-kriging & 98.01\% & 0.038 & 0.14  & 0.046  \\
\hline
\end{tabular}
\end{center}
\caption{Application: hydrodynamic simulator. Comparison between 
kriging and co-kriging. The co-kriging model provides predictions significantly better than the ones of the kriging model.}
\label{tab4}
\end{table}
We can see that the difference of accuracy between the two models is important. Indeed, the one of the co-kriging model is significantly better. Furthermore, comparing the RMSE and the RIMSE estimations  in Table \ref{tab4}, we see that we have a good estimation of the predictive distribution variances for the two models. We note that the predictive variance for the co-kriging is obtained with a simple co-kriging model. Therefore, it will be slightly larger in the universal co-kriging case. Indeed, by computing the universal co-kriging equations, we find $\mathrm{RIMSE} = 0.058$.

We can compare the RMSE obtained with the test set with the RMSE obtained with a Leave-One-Out cross validation procedure (see Section \ref{FastCV}). For this procedure, we test our model on $n_2 = 5$ validation sets obtained by removing one observation at a time. As presented in Section  \ref{FastCV}, we can either choose to remove the observations from $z_2$ or from $z_2$ and $z_1$. The root mean squared error of the  Leave-One-Out cross validation procedure obtained by removing observations from $z_2$ is RMSE$_{z_2,LOO} = 4.80.10^{-3}$ whereas the one obtained by removing observations from $z_2$ and $z_1$ is RMSE$_{z_1,z_2,LOO} = 0.10$. Comparing RMSE$_{z_2,LOO} $ and RMSE$_{z_1,z_2,LOO}$ to the RMSE obtained with the external test set, we see that the procedure which consists in removing points from $z_2$ and $z_1$ provides a better proxy for the generalization error. Indeed, RMSE$_{z_2,LOO} $ is a relevant proxy for the generalization error only at points where $z_1$ is available. Therefore, it underestimates the error at locations where $z_1$ is unknown.

\begin{figure}[h]
\begin{center}
\includegraphics[width =7cm]{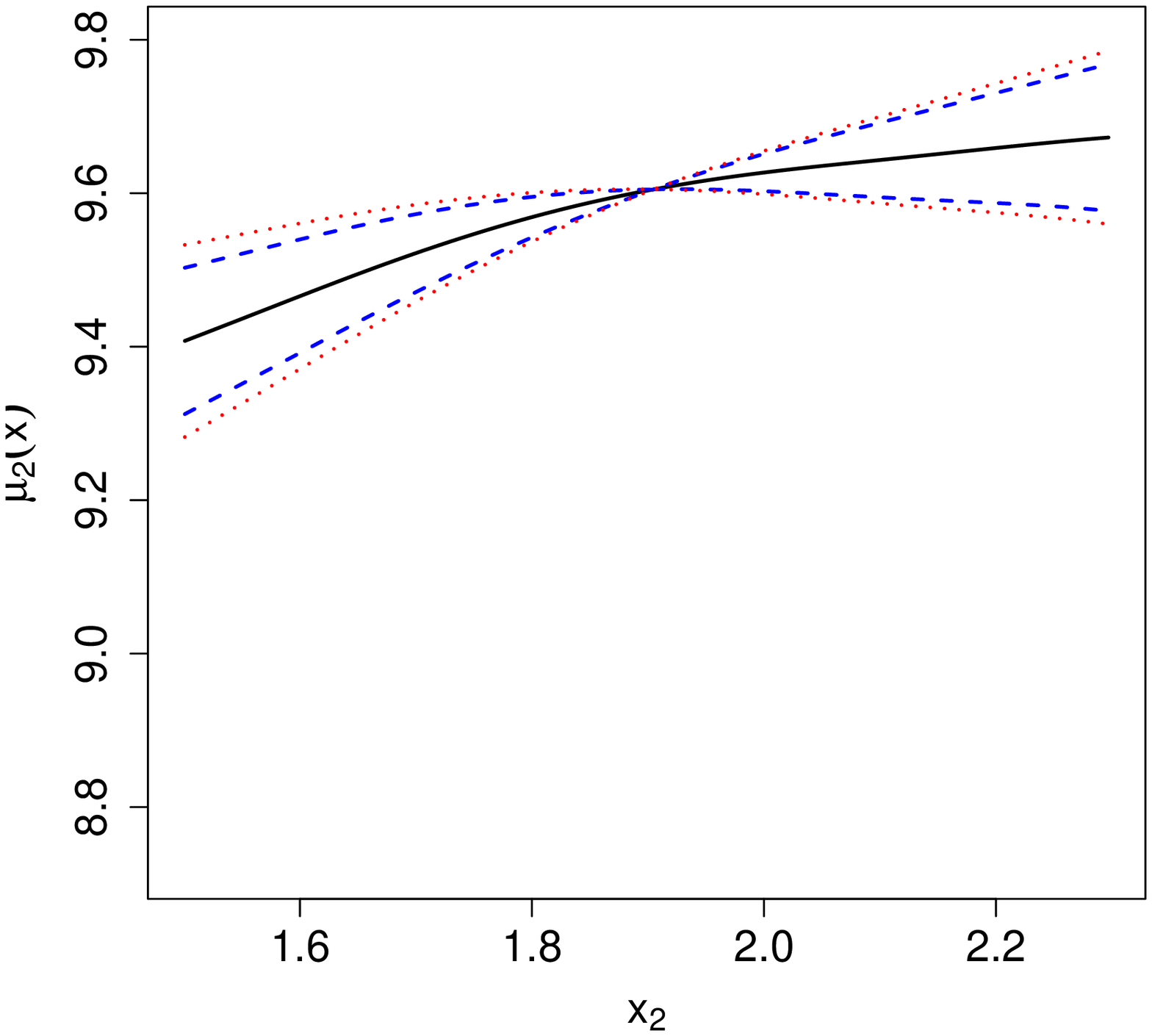}
\includegraphics[width = 7cm]{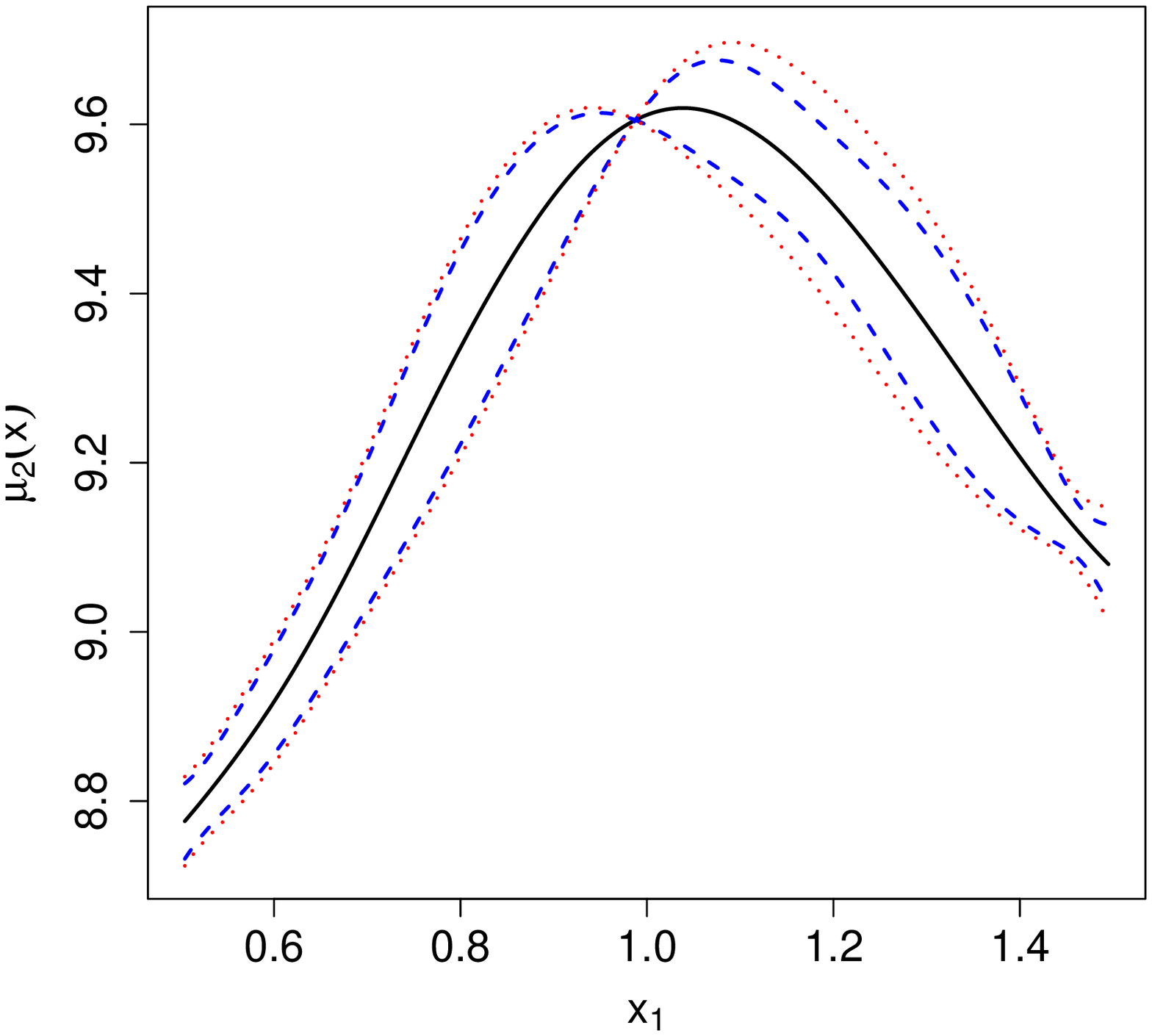}
\caption{Mean and confidence intervals for the simple and the universal co-kriging. The figure on the left hand side represents the predictions along the vertical line $x_1=0.99$ and the figure on the right hand side represents the predictions along the horizontal line $x_2=1.91$. The solid black lines represent the mean of the two co-kriging models, the dashed lines represent the confidence interval at plus or minus twice the standard deviation of the simple co-kriging and the dotted lines represent the same confidence intervals for the universal co-kriging.}
\label{fig:predline}
\end{center}
\end{figure}

Figure \ref{fig:predline} represents the mean and  confidence intervals at plus or minus twice the standard deviation of the simple and universal co-krigings for points along the vertical line $x_1 = 0.99$ and the horizontal line $x_2=1.91$ ($x=(0.99,1.91)$ corresponds to the coordinates of the point  of $D_2$ in the center of the domain $ [0.5,1.5]\times[1.5,2.3]$ in Figure \ref{fig:plan}). In Figure \ref{fig:predline} on the right hand side, we see a  necked point around the coordinates $x_1 = 1.5$  since, in the direction of $x_2$, the  correlation hyper-parameters  length for $Z_1(x)$ and $ \delta_2(x) $ are large ($\theta_{1,2} = 1.20$ and $\theta_{2,2} = 1.37$) and a  point of $D_2$ has almost the same coordinate.

\subsection{Multi-fidelity surrogate model for the mixture characteristic length $L_c$}

In this section, we build a co-kriging model for the mixture characteristic length $L_c$. The aim of this example is to highlight that it can  be worth having an adjustment coefficient $\rho_1$ depending on $x$. We use the same training and test sets as in the previous section and we consider a tensorised matern-$\frac{5}{2}$ kernel (\ref{mat2D}). Let us consider the two following cases:
\begin{itemize}
\item{Case 1:} $g_1(x) = 1$, $f_2(x) = 1$ and $f_1(x)=1$
\item{Case 2:} $g_1^T(x) = (\begin{array}{cc}1 &  x_1\end{array})$, $f_2(x) = 1$ and $f_1(x)=1$
\end{itemize}
We have the following hyper-parameter maximum likelihood estimates for the two cases
\begin{itemize}
\item{Case 1:}  $\hat{\theta}_1 = (0.52, 1.09) $ and  $\hat{\theta}_2 = (0.03, 0.02)$
\item{Case 2:}  $\hat{\theta}_1 = (0.52, 1.09) $ and  $\hat{\theta}_2 = (0.14, 1.37)$
\end{itemize}
The estimation of $\hat{\theta}_1$ is identical in the two cases since it does not depend on $\rho_1$ and it is estimated with the same observations. Furthermore, we can see an important difference between the estimates of $\hat{\theta}_2 $. Indeed, they are larger in the Case 2 than in the Case 1 which indicates  that the model is smoother in the Case 2.
Table \ref{tab1Lc} presents the estimations of $\beta_1$ and $\sigma_1^2$ for the two cases (see Section \ref{paramestim}).

\begin{table}[H]
\begin{center}
\begin{tabular}{|c|c|c|}
\hline
Trend coefficient & $\Sigma_1 \nu_1$ & $\Sigma_1/\sigma_1^2$ \\
\hline
$\beta_1$ & $1.26$ & $ 0.97$ \\
\hline
\hline
Variance coefficient & $Q_1$ & $2 \alpha_1$ \\
\hline
$\sigma_1^2$ & $15.62$ & $24$ \\
\hline
\end{tabular}
\end{center}
\caption{Application: hydrodynamic simulator. Estimations of $\beta_1$ and $\sigma_1^2$  for the response $L_c$ (see equations (\ref{eq8}) and (\ref{eq9})).}
\label{tab1Lc}
\end{table}
Then, Table \ref{tab2Lc} presents the estimations of $\beta_2$, $ \beta_{\rho_1}$ and $\sigma_2^2$ for the Case 1, i.e. when $\rho_1$ is constant  (see Section \ref{paramestim}).

\begin{table}[H]
\begin{center}
\begin{tabular}{|c|c|c|}
\hline
Trend coefficient & $\Sigma_2 \nu_2$ & $\Sigma_2/\sigma_2^2$ \\
\hline
$\left( \begin{array}{c} \beta_{\rho_1} \\ \beta_2 \\ \end{array} \right)$ & $ \left( \begin{array}{c} 1.49  \\ -0.26  \\ \end{array} \right)$ & $ \left( \begin{array}{cc}  0.83   & -0.79  \\  -0.79 & 0.95 \\  \end{array} \right) $\\
\hline
\hline
Variance coefficient & $Q_2$ & $2 \alpha_2$ \\
\hline
$\sigma_2^2$ & $0.01$ & $3$ \\
\hline
\end{tabular}
\end{center}
\caption{Application: hydrodynamic simulator.  Estimations of  $\beta_2$, $ \beta_{\rho_1}$ and $\sigma_2^2$ for the Case 1, i.e. when  $\rho_1$ is  constant, for the response $L_c$     (see equations (\ref{eq8}) and (\ref{eq9})).}
\label{tab2Lc}
\end{table}

Finally, Table \ref{tab3Lc} presents the estimations of $\beta_2$, $ \beta_{\rho_1}$ and $\sigma_2^2$ for the Case  2, i.e. when $\rho_1$  depends on $x$  (see Section \ref{paramestim}).

\begin{table}[H]
\begin{center}
\begin{tabular}{|c|c|c|}
\hline
Trend coefficient & $\Sigma_2 \nu_2$ & $\Sigma_2/\sigma_2^2$ \\
\hline
$\left( \begin{array}{c} \beta_{\rho_1} \\ \beta_2 \\ \end{array} \right)$ & $ \left( \begin{array}{c} 1.66 \\ -0.48  \\ -0.04  \\ \end{array} \right)$ & $ \left( \begin{array}{ccc}  2.34   & -3.50  & 0.44  \\  -3.50 &  9.18 &   -3.67  \\ 0.44 & -3.67 &   2.60  \end{array} \right) $\\
\hline
\hline
Variance coefficient & $Q_2$ & $2 \alpha_2$ \\
\hline
$\sigma_2^2$ & $3.24.10^{-4}$ & $2$ \\
\hline
\end{tabular}
\end{center}
\caption{Application: hydrodynamic simulator.  Estimations of  $\beta_2$, $ \beta_{\rho_1}$ and $\sigma_2^2$ for the Case 2, i.e. when  $\rho_1$   depends on $x$, for the response $L_c$     (see equations (\ref{eq8}) and (\ref{eq9})).}
\label{tab3Lc}
\end{table}
We see in Table \ref{tab2Lc}   that the adjustment coefficient is around 1.5 which indicates  that the magnitude of the expensive code is slightly more important than the one of the cheap code. Furthermore, we see in Table \ref{tab3Lc} that if we consider an adjustment coefficient which linearly depends on $x_1$ (i.e. with $g_1^T(x) =  (\begin{array}{cc}1 &  x_1\end{array})$), the constant part of $\beta_{\rho_1}$  is  more important (it is around 1.66) and there is a negative slope in the direction $x_1$ (it is around $-0.48$). Since $x \in [0.5, 1.5]$, the averaged value of $\rho_1$ is 1.18 and goes from 1.42 at $x_1 = 0.5$ to 0.94 at $x_1 = 1.5$.  We see also a significant difference between the two case for the variance estimation. Indeed, the variance estimate in the Case 1 (see Table \ref{tab2Lc}) is much more important than the one in the Case 2 (see Table \ref{tab3Lc}). This could mean that we learn better in the Case 2 than in the Case 1.

\begin{figure}[h]
\begin{center}
\includegraphics[width =7cm]{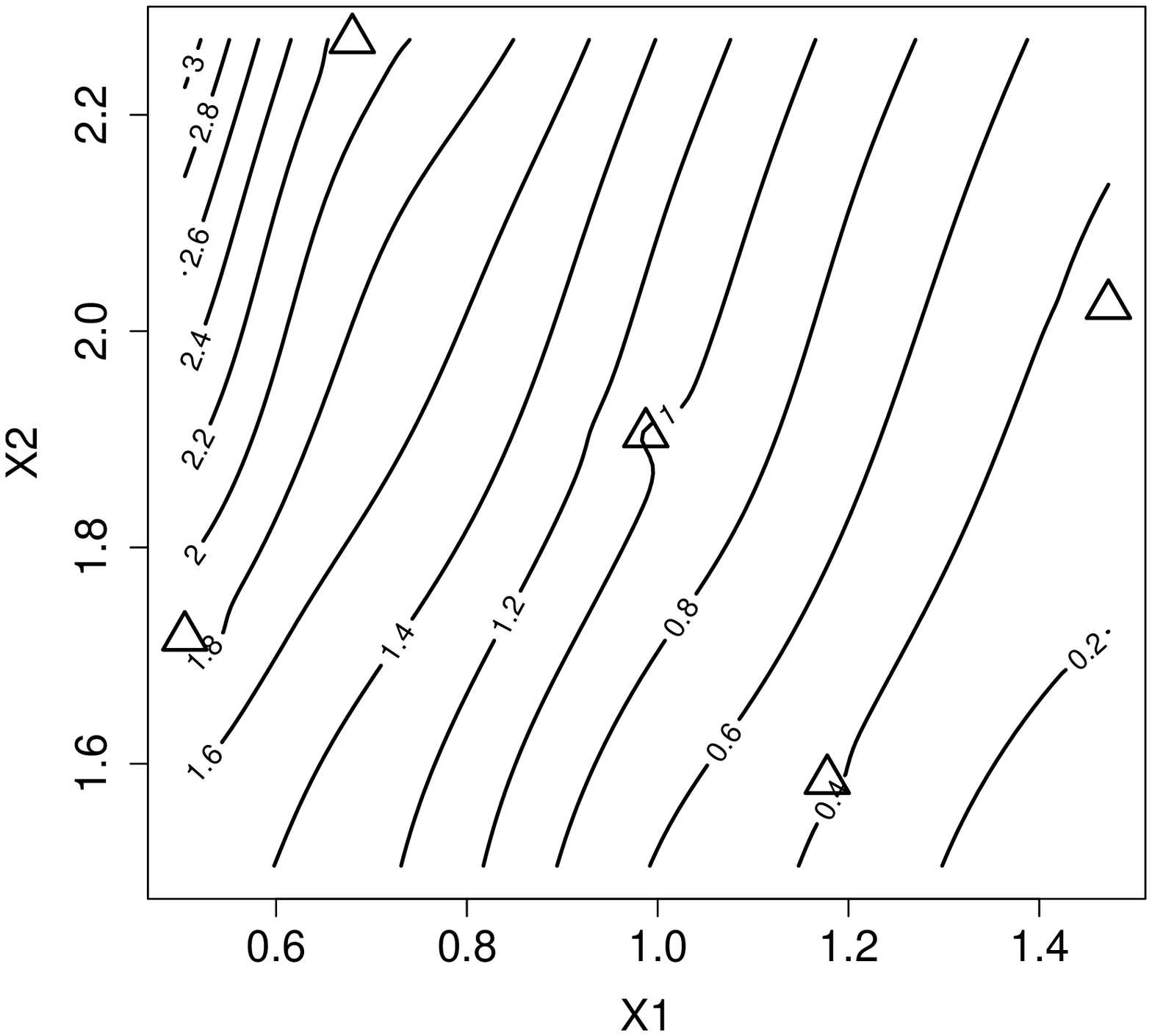}
\includegraphics[width = 7cm]{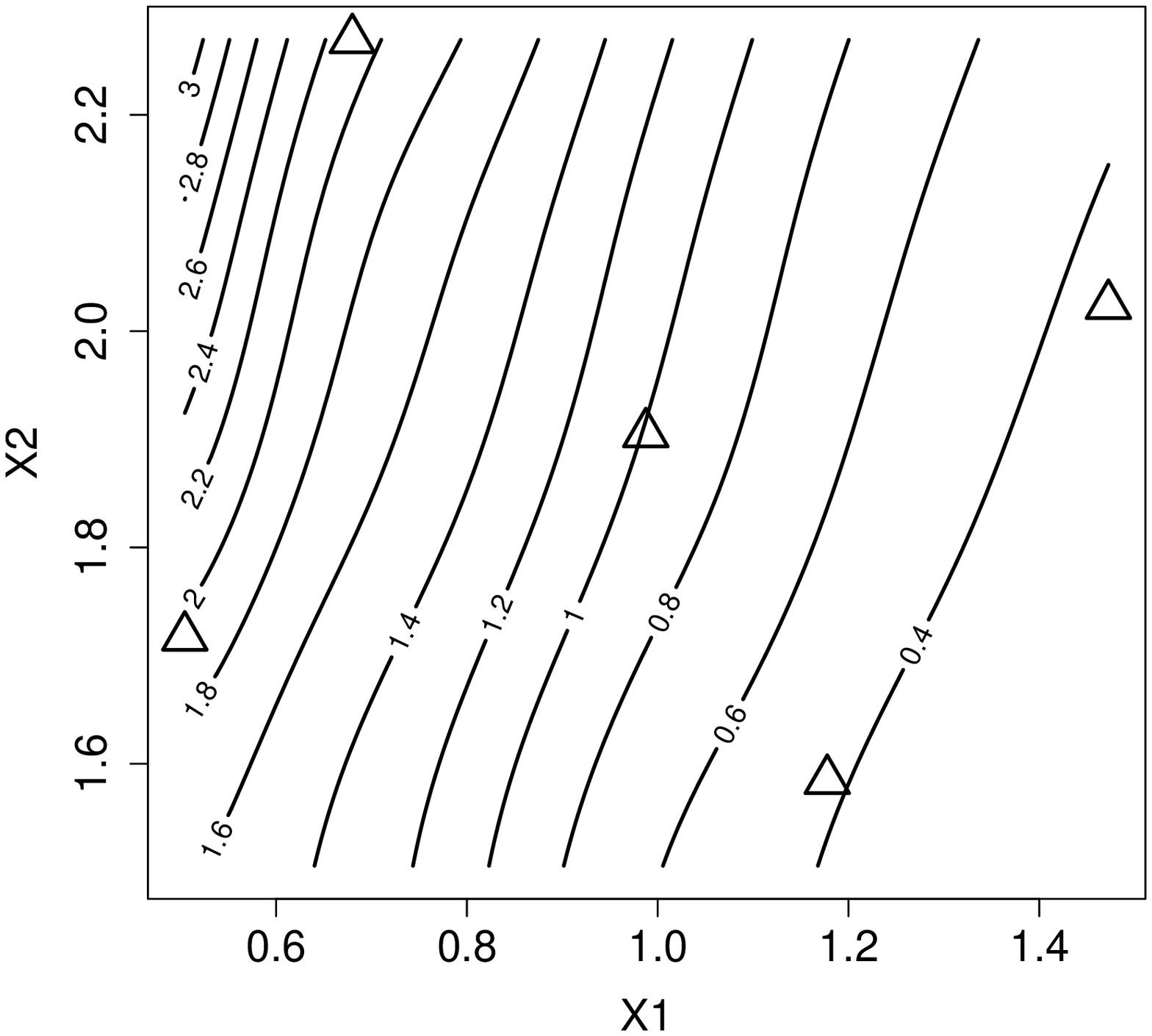}
\caption{Contour plot of the co-kriging mean when $\rho_1$ is constant (on the left hand side) and  when $\rho_1$ is depends on $x$ (of the right hand side). The triangles represent the $n_2 = 5$ points of the experimental design set  of the expensive code.}
\label{fig:contourLc}
\end{center}
\end{figure}

Figure \ref{fig:contourLc} illustrates the contour plot of the two co-kriging models, i.e. when $\rho_1$ is constant and when $\rho_1$ depends on $x$. 

Furthermore, Table \ref{tab5} compares the prediction accuracy of the co-kriging in the two cases. The precision is computed on the test set of 175 observations.

\begin{table}[H]
\begin{center}
\begin{tabular}{ccc}
\hline
& RMSE & MaxAE    \\
Case 1 & $7.26.10^{-3}$  & 0.23  \\
Case 2 &  $1.53.10^{-3}$ & 0.16  \\
\hline
\end{tabular}
\end{center}
\caption{Application: hydrodynamic simulator. Comparison between 
co-kriging when $\rho_1$ is constant (Case 1) and co-kriging when $\rho_1$ depends on $x$ (Case 2). The Case 2 provides predictions  better than the Case 1, it is hence worthwhile to consider an adjustment coefficient that is not constant.}
\label{tab5}
\end{table}

We see that the co-kriging model in Case 2 is clearly  better than the one in  Case 1. Therefore, we illustrate in this application that  it can be  worth   considering an adjustment coefficient not constant contrarily to the model presented in \cite{KO00} and \cite{AF07}.

\section{Conclusion}

We have  presented in this paper a recursive formulation for a multi-fidelity co-kriging model. This model allows us to build  surrogate models using data from simulations  of different levels of fidelity. 

The strength of the suggested approach is that it  considerably reduces  the complexity    of the  co-kriging model    while it preserves its  predictive efficiency.  Furthermore, one of the most important consequences of the recursive formulation is that the construction of the surrogate model is equivalent to build $s$ independent krigings. Consequently, we can naturally adapt results of kriging to the proposed co-kriging model.

First, we present a Bayesian estimation of the model   parameters  which provides closed form expressions  for the parameters of the posterior distributions. We note that, from these posterior distributions, we can deduce the maximum likelihood estimates of the parameters. Second, thanks to the joint distributions  of the parameters and the recursive formulation, we can deduce closed form formulas for the mean and covariance of the posterior predictive distribution. Due to their similarities with the universal kriging equations, we call these formulas the universal co-kriging equations.  Third, we present closed form expressions for  the cross-validation equations of the co-kriging surrogate model. These expressions reduce considerably the complexity of the cross-validation procedure and are derived from the one of kriging model that we have extended. 

The suggested model has been successfully applied  to a hydrodynamic code. We also present in this application a practical way to design the experiments of the multi-fidelity model.

\section{Acknowledgements}

The authors   thank  Dr. Claire Cannamela   for providing the data for the application and for     interesting discussions.

\appendix

\section{Proofs}

\subsection{Proof of Proposition 1}\label{A1}

Let us consider the co-kriging mean of the model (\ref{eq1})  presented in \cite{KO00} for a $t$-level co-kriging with $t = 2,\dots,s$:

\begin{displaymath}
m_{Z_t}(x) = h_t(x)^T\beta^{(t)} + t_t(x)^TV_t^{-1}(z^{(t)}-H_t\beta^{(t)})
\end{displaymath}
where  $\beta^{(t)} = (\beta_1^T,\dots,\beta_t^T)^T$, $z^{(t)} = (z_1^T,\dots,z_t^T)^T$ and $h_t(x)^T$ is defined in the following equation:
\begin{equation} 
h_t(x)^T = \left( \left(\prod_{i=1}^{t-1} {\rho_{i}(x) }\right)f_1^T(x) , \left(\prod_{i=2}^{t-1} {\rho_{i}(x) }\right)f_2^T(x) , \dots , \rho_{t-1}(x) f_{t-1}^T(x) , f_t^T(x) \right)
\end{equation}
We have:
\begin{eqnarray*}
h_t(x)^T\beta^{(t)}  & = &  \rho_{t-1}(x) \left(
\left(\prod_{i=1}^{t-2} {\rho_{i}(x) }\right)f_1^T(x) , \left(\prod_{i=2}^{t-2} {\rho_{i}(x) }\right)f_2^T(x) , \dots ,  f_{t-1}^T(x)
\right)\beta^{(t-1)}
+ f_t^T(x) \beta_t \\
 & = & \rho_{t-1}(x)h_{t-1}(x)^T\beta^{(t-1)} + f_t^T(x) \beta_t  
\end{eqnarray*}
Then, from equations:
\begin{equation} 
\mathrm{cov}(Z_t(x),Z_ {t'}(x')|\sigma^2,\beta,\beta_\rho) =
\left(\prod_{i=t'}^{t-1}\rho_i(x)\right) \mathrm{cov}(Z_{t'}(x),Z_{t'}(x')|\sigma^2,\beta,\beta_\rho)
\end{equation}
and:
\begin{equation} 
\mathrm{cov}(Z_t(x),Z_t(x')|\sigma^2,\beta,\beta_\rho) =
\sum_{j=1}^{t}{\sigma_{j}^2\left( \prod_{i=j}^{t-1} {\rho_{i}(x) \rho_{i}(x')} \right)r_j(x,x')}
\end{equation}
with $t > t'$ , we have the following equality: 
\begin{eqnarray*}
t_t(x)^TV_t^{-1}z^{(t)}& = & \rho_{t-1}(x) t_{t-1}(x)^TV_{t-1}^{-1}z^{(t-1)} -\left(\rho_{t-1}^T(D_t)\right)\odot \left(r_t^T(x)R_t^{-1}z_{t-1}(D_t)\right) \\
 &  &+ r_t^T(x) R_t^{-1}z_t 
\end{eqnarray*}
and with equation (\ref{eq6}):
\begin{displaymath}
t_t(x)^TV_t^{-1}H_t\beta^{(t)} = \rho_{t-1}(x) t_{t-1}(x)^TV_{t-1}^{-1}H_{t-1}\beta^{(t-1)}+
 r_t^T(x)R_t^{-1} F_t(D_t) \beta_t 
\end{displaymath}
where $\odot$ stands for the element by element matrix product.
We hence obtain the recursive relation:
\begin{displaymath}
m_{Z_t}(x)  =\rho_{t-1}(x) m_{Z_{t-1}}(x) + f_t^T(x) \beta_t +  r_t^T(x)R_t^{-1}\left[
z_t  - \rho_{t-1}(D_t)\odot  z_{t-1}(D_t) -  F_t(D_t) \beta_t 
\right]
\end{displaymath}
The co-kriging mean of the model (\ref{eq7ter}) satisfies the same recursive relation (\ref{eq6}), and we have  $m_{Z_1}(x) = \mu_{Z_1}(x)$. This proves the first equality of Proposition 1:
\begin{displaymath}
\mu_{Z_{s}}(x)  =    m_{Z_{s}}(x)
\end{displaymath}

We follow the same guideline for the co-kriging covariance:
\begin{displaymath}
s^2_{Z_t}(x,x') = v^2_{Z_t}(x,x') - t_t^T(x) V_t^{-1}t_t(x')
\end{displaymath}
where $v^2_{Z_t}(x,x')$ is the covariance between $Z_t(x)$ and $Z_t(x')$ and $s^2_{Z_t}(x,x') $ is  the covariance function of the conditioned Gaussian process $[Z_t(x)|\mathcal{Z}^{(t)}=z^{(t)},\beta,\beta_{\rho},\sigma^2]$ for the model (\ref{eq1}).
 From equation (\ref{eq7bis}), we can deduce the following equality:
\begin{displaymath}
 \sigma^2_{Z_t}(x,x')  = \rho_{t-1}(x)\rho_{t-1}(x') v^2_{Z_{t-1}}(x,x') +v_t^2 (x,x')
\end{displaymath}
where $\sigma^2_{Z_t}(x,x')$ is the covariance function of the conditioned Gaussian process $[Z_t(x)|\mathcal{Z}^{(t)}=z^{(t)},\beta_t,\beta_{\rho_{t-1}},\sigma_t^2]$ of the recursive model (\ref{eq7ter}).
 Then, from equation (\ref{eq7}) and (\ref{eq7bis}), we have:
\begin{displaymath}
t_t^T(x) V_t^{-1}t_t(x') = \rho_{t-1}(x)\rho_{t-1}(x')t_{t-1}^T(x) V_{t-1}^{-1}t_{t-1}(x')+\sigma_t^2  r_t^T(x)R_t^{-1}r_t(x') 
\end{displaymath}
Finally we can deduce the following equality:
\begin{displaymath}
s^2_{Z_t}(x,x') =   \rho_{t-1}(x)\rho_{t-1}(x')\left(v^2_{Z_{t-1}}(x,x')  - t_{t-1}^T(x) V_{t-1}^{-1}t_{t-1}(x')\right)+\sigma_t^2\left(1-r_t^T(x)R_t^{-1}r_t(x')\right) 
\end{displaymath}
which is equivalent to:
\begin{displaymath}
s^2_{Z_t}(x,x')= \rho_{t-1}(x)\rho_{t-1}(x')s^2_{Z_{t-1}}(x,x')+\sigma_t^2\left(1-r_t^T(x)R_t^{-1}r_t(x')\right) 
\end{displaymath}
This is the same recursive relation as the one satisfies by the co-kriging covariance $\sigma_{Z_t}^2(x,x')$ of the model  (\ref{eq7ter}) (see equation (\ref{eq12})). Since $s^2_{Z_1}(x,x') = \sigma^2_{Z_1}(x,x')$, we have :
\begin{displaymath}
\sigma_{Z_{s}}^2(x,x')   =   s_{Z_{s}}^2(x,x') 
\end{displaymath}
This  equality with $x=x'$ proves  the second equality of Proposition 1. \hfill $\Box$

\subsection{Proof of Proposition 2}\label{A2}

Noting that the mean of the predictive distribution in equation (\ref{eq10}) does  not depend on $\sigma_t^2$ and thanks to the law of total expectation, we have the following  equality: 
\begin{displaymath}
\mathbb{E} \left[Z_{t}(x)|\mathcal{Z}^{(t)}=z^{(t)} \right]  = \mathbb{E} \left[ \mathbb{E} \left[Z_{t}(x)|\mathcal{Z}^{(t)}=z^{(t)},\sigma^2_t, \beta_t,\beta_{\rho_{t-1}}\right]\left|\mathcal{Z}^{(t)}=z^{(t)} \right.  \right] 
\end{displaymath}
From the equations (\ref{eq11}) and  (\ref{eq8}), we directly deduce the equation (\ref{13}). Then, we have the following equality:
\begin{equation}\label{covbrho}
\mathrm{var} \left(
\mu_{Z_t}(x) \left|z^{(t)} ,\sigma_t^2 \right. \right) = (h_t^T(x)-r_t(x)^TR_t^{-1}H_t)\Sigma_t (h_t^T(x)-r_t(x)^TR_t^{-1}H_t)^T
\end{equation}
The law of total variance states that:
\begin{eqnarray*}
 \mathrm{var}(Z_t(x) |z^{(t)} ,\sigma_t^2) & = &  \mathbb{E}  \left[ \mathrm{var}(Z_t(x) |z^{(t)} ,\beta_t,\beta_{\rho_{t-1}},\sigma_t^2) \left|z^{(t)} ,\sigma_t^2 \right. \right]  \\
&+ &\mathrm{var}\left(
\mathbb{E} \left[ Z_t(x) |z^{(t)} ,\beta_t,\beta_{\rho_{t-1}},\sigma_t^2 \right]
\left|z^{(t)} ,\sigma_t^2 \right. \right)
\end{eqnarray*}
Thus, from equations (\ref{eq11}),  (\ref{13}) and  (\ref{covbrho}), we obtain:
\begin{equation}
\begin{array}{lll}
\mathrm{var}(Z_t(x)|\mathcal{Z}^{(t)} = z^{(t)},\sigma_t^2) &  = &  \hat{\rho}_t^2(x) \mathrm{var}(Z_{t-1}(x )|\mathcal{Z}^{(t-1)} = z^{(t-1)},\sigma_t^2 )+ \sigma_t^2 \left(1-r_t^T(x)R_t^{-1}r_t^T(x)\right) \\
 & &+ \left(
h_t^T - r_t^T(x)R_t^{-1}H_t\right) \Sigma_t \left(
h_t^T - r_t^T(x)R_t^{-1}H_t\right)^T \\
\end{array}
\end{equation}
Again using the law  of total variance and the independence between $\mathbb{E} \left[Z_{t}(x)|\mathcal{Z}^{(t)}=z^{(t)}, \beta_t,\beta_{\rho_{t-1}}\right] $ and $\sigma_t^2$, we have:
\begin{equation} 
 \mathrm{var}(Z_t(x) |z^{(t)})  = \mathbb{E} \left[  \mathrm{var}(Z_t(x)) |z^{(t)},\sigma_t^2\right]
\end{equation}
We obtain the equation (\ref{19}) from equation (\ref{eq9}) by noting that the mean of an inverse Gamma distribution $\mathcal{IG}(a,b)$ is $b/(a-1)$. \hfill $\Box$

\subsection{Proof of Proposition 3}\label{A3}

Let us consider that  $\xi_s$ is the index of the $k$ last points of $D_s$. We denote  by $D_\mathrm{test}$ these points. 
First we consider the variance and the trend parameters as fixed, i.e. $\sigma_{t,-\xi_t}^2 = \frac{Q_t}{2(a_t-1)}$ and $\lambda_{t,-\xi_t} = \Sigma_t \nu_t$, and $\mathcal{V}_s = 0$, i.e. we are in  the simple co-kriging case.
Thanks to the block-wise inversion formula, we have the following equality:
\begin{equation}\label{blockinverse}
R_s^{-1} = \left(
\begin{array}{cc}
A & B \\
B^T & \mathcal{Q}^{-1} \\
\end{array}
\right)
\end{equation}
with $A = \left[ R_s^{-1} \right]_{[-\xi_s,-\xi_s]} + \left[ R_s^{-1} \right]_{[-\xi_s,-\xi_s]} \left[ R_s^{-1} \right]_{[-\xi_s,\xi_s]}Q^{-1}\left[ R_s^{-1} \right]_{[\xi_s,-\xi_s]}\left[ R_s^{-1} \right]_{[-\xi_s,-\xi_s]} $,\\  $B=-\left[ R_s^{-1} \right]_{[-\xi_s,-\xi_s]} \left[ R_s^{-1} \right]_{[-\xi_s,\xi_s]}\mathcal{Q}^{-1}$ and:
\begin{equation}\label{Q}
 \mathcal{Q} =  \left[ R_s^{-1} \right]_{[\xi_s,\xi_s]} - \left[ R_s^{-1} \right]_{[ \xi_s, - \xi_s]} \left(\left[ R_s^{-1} \right]_{[-\xi_s,-\xi_s]} \right)^{-1} \left[ R_s^{-1} \right]_{[-\xi_s, \xi_s]}
\end{equation}
We note that  $\frac{Q_s}{2(a_s-1)} \mathcal{Q} =  \frac{Q_t}{2(a_t-1)}  \left(\left[R_s^{-1}\right]_{[\xi_s,\xi_s]}\right)^{-1}$ represents the covariance  matrix of the  points in $D_\mathrm{test}$ with respect to the covariance kernel of  a Gaussian process of kernel $\frac{Q_s}{2(a_s-1)} r_s(x,x')$ (which is the one of $\delta_s(x)$) conditioned by the points $D_s \setminus D_\mathrm{test}$.  Therefore, from the previous remark and the equation (\ref{eq12}), we can deduce the equation (\ref{CVvar}).

Furthermore, we have the following equality:
\begin{equation}
\begin{array}{lll}
\left( \left[ R_s^{-1} \right]_{[\xi_s,\xi_s]} \right)^{-1} \left[ R_s^{-1}\left(
z_s - H_s \lambda_{s,-\xi_s}
\right) \right]_{[\xi_s ]} & =& 
z_s(D_\mathrm{test}) - h_s^T(D_\mathrm{test}) \Sigma_s \nu_s \\
 & - & 
  \left[ R_s^{-1} \right]_{[-\xi_s,\xi_s]}\left(   \left[ R_s^{-1} \right]_{[\xi_s,\xi_s]} \right)^{-1} \\
 & \times & \left( z_{s-1}(D_s \setminus D_\mathrm{test}) - [ H_s^T]_{[- \xi_s]}\Sigma_s \nu_s  \right)
\end{array}
\end{equation}
From this equation and equation (\ref{eq11}), we can directly deduce the equation (\ref{CVerr}) with $\varepsilon_{Z_s,\xi_s} =  z_s(D_\mathrm{test}) - \mu_{Z_s}(D_\mathrm{test}) $. 

Then, we suppose the trend and the variance parameters as unknown and we have to re-estimate them when we remove the observations. Thanks to the parameter estimations presented in Section \ref{paramestim}, we can deduce that the estimates of $\sigma_{t,-\xi_t}^2 $ and $\lambda_{t,-\xi_t} $ when we remove observations of index $\xi_t$ are given by the following equations:
\begin{equation}
\lambda_{s,-\xi_s} \left([H_s^T ]_{-\xi_s}K_s [H_s]_{-\xi_s} \right)= [H_s^T]_{-\xi_s} K_s z_s(D_s \setminus D_{test})
\end{equation}
and:
\begin{equation} 
 \sigma_{s,-\xi_s}^2 = \frac{\left(z_s(D_s \setminus D_{\mathrm{test}}) - [H_s]_{-\xi_s}\lambda_{s,-\xi_s}  \right)^T K_s \left(z_s(D_s \setminus D_{\mathrm{test}}) - [H_s]_{-\xi_s}\lambda_{s,-\xi_s}  \right)  }{n_s-p_s-q_{s-1}-n_{train}}
\end{equation}
with $K_s =  \left( \left[ R_s \right]_{[-\xi_s,-\xi_s]} \right)^{-1}$. 

From the equality (\ref{blockinverse}), we can deduce that $K_s = A - B \mathcal{Q} B^T$ from which we  obtain the equation (\ref{KtCV}).
Finally, to obtain the cross-validation equations for the universal co-kriging, we just have to estimate the following quantity (see equation (\ref{19})):
\begin{equation} 
\left(h_s^T(D_\mathrm{test})^T -  \left[ R_s^{-1} \right]_{[-\xi_s,\xi_s]}K_s[H_s]_{-\xi_s}\right) \Sigma_s \left(
h_s^T(D_\mathrm{test})^T -   \left[ R_s^{-1} \right]_{[-\xi_s,\xi_s]}K_s[H_s]_{-\xi_s}\right)^T 
\end{equation}
with $\Sigma_s = \left([H_s^T ]_{-\xi_s}K_s [H_s]_{-\xi_s} \right)^{-1}$. The following equality:
\begin{equation} 
\left(h_s^T(D_\mathrm{test})^T -  \left[ R_s^{-1} \right]_{[-\xi_s,\xi_s]}K_s[H_s]_{-\xi_s}\right) =   \left( \left( [R_s^{-1}]_{[\xi_s,\xi_s]}\right)^{-1}
\left[ R_s^{-1}H_s \right]_{[\xi_s ]} 
\right)
\end{equation}
allows us to obtain the equation (\ref{CVUnivCok}) and completes the proof. \hfill $\Box$

\bibliographystyle{apalike}
\bibliography{biblio}

\label{fin}

\end{document}